\documentclass[12pt]{article}
\usepackage{amssymb}
\usepackage{amsmath}
\usepackage{indentfirst}

\newcommand{\ba}{\begin{array}}
\newcommand{\ea}{\end{array}}
\newcommand{\bi}{\begin{itemize}}
\newcommand{\ei}{\end{itemize}}
\newcommand{\bc}{\begin{center}}
\newcommand{\ec}{\end{center}}
\newcommand{\bfr}{\begin{flushright}}
\newcommand{\efr}{\end{flushright}}
\newcommand{\f}{\frac}
\newcommand{\ov}{\overline}

\newcommand{\ds}{\displaystyle}

\newcommand{\q}{\quad}

\begin{document}

\title{TRIGONOMETRIC AND HYPERBOLIC INEQUALITIES}
\author{J\'ozsef S\'andor\\
Babe\c{s}-Bolyai University\\
Department of Mathematics\\
Str. Kog\u{a}lniceanu nr. 1\\
400084 Cluj-Napoca, Romania\\
email: jsandor@math.ubbcluj.ro}
\date{}
\maketitle

\begin{abstract}
This is a collection of 10 papers of the author  on various trigonometric and hyperbolic inequalities.
One can find here Kober and Jordan type  trigonometric  inequalities,
Wilker, Huygens, Cusa-Huygens type inequalities,  as well as trigonometric-hyperbolic
inequalities, connecting trigonometric and hyperbolic functions in the same inequality.
Best possible results are pointed out, too.
A common point in many arguments  is the use of the bivariate means and their inequalities,
or the monotonicity and convexity properties of certain auxiliary functions.
\end{abstract}


\newpage
\setcounter{section}{0}

\bc
{\Large\bf 1. On New Refinements of Kober's and Jordan's Trigonometric Inequalities}
\ec

\begin{abstract}
This paper deals with some inequalities for trigonometric functions such as the Jordan
inequality and Kober's inequality with refinements.
In particular, lower and upper bounds for functions such as
$(\sin x)/x$, $(1-\cos x)/x$ and $(\tan x/2)/x$ are proved.
\end{abstract}



\section{Introduction}

During the past several years there has been a great interest in trigonometric inequalities
(see e.g. a list of papers in \cite{54-2}, where 46 items are included).

The classical Jordan inequality
$$\ds\f{2}{\pi }x\le \sin x,\q
0\le x\le \ds\f{\pi }{2}
\eqno(1)$$
has been in the focus of these studies and many refinements have been proved
(see e.g. the references from \cite{54-2}, \cite{54-12}, \cite{54-13}, \cite{54-4}).
One of the early references from this topic is the author's book \cite{54-7}
on Geometric inequalities (with a geometric interpretation of (1), too, which was republished
also in the book \cite{54-16}) or the author's papers from 2001 (\cite{54-9}) and 2005 (\cite{54-10}, \cite{54-11}),
or more recently, from 2007 (\cite{54-13}).
See also the book \cite{54-15}.
For example, in \cite{54-9} it was proved that
$$\ds\f{1+\cos x}{2}<\ds\f{\sin x}{x}<\cos \ds\f{x}{2},\q
0<x<\ds\f{\pi }{2},
\eqno(2)$$
rediscovered, and used many times in the literature (see e.g. \cite{54-2}, relations (1.2) and (1.10)).

Another famous inequality, connected with Jordan's inequality is Kober's inequality
(see e.g. \cite{54-3})
$$\cos x\ge 1-\ds\f{2}{\pi }x,\q 0\le x\le \ds\f{\pi }{2}.
\eqno(3)$$

Though (3) can be proved by considering the monotonic property of the function
$\ds\f{1-\cos x}{x}$,
in \cite{54-10}, \cite{54-11}, \cite{54-13} we have remarked that, it follows at one via the substitution
$x\to \ds\f{\pi }{2}-x$
in relation (1);
and vice-versa, (3) implies (1) by the same manner.
In paper \cite{54-11} we have shown that the application
$g(x)=\sin x/x$, $0<x\le \ds\f{\pi }{2}$, $g(0)=1$,
is a strictly increasing and strictly concave function on $\left[0,\ds\f{\pi }{2}\right]$.
See also \cite{54-14}.
By writing that the graph of line passing on the points $A(0,1)$ and
$B\left(\ds\f{\pi }{2},\ds\f{2}{\pi }\right)$ from the graph of $g$ is below
the graph of $g$, one obtains that
$$\ds\f{\sin x}{x}\ge 1+\ds\f{2(2-\pi )}{\pi ^2}\cdot x,
\eqno(4)$$
which is a refinement of (1), as the right side of (4) may be written also as
$$\ds\f{2}{\pi }+\ds\f{\pi -2}{\pi ^2}(\pi -2x)\ge \ds\f{2}{\pi }.$$

By writing the tangent line to the graph of function $g$ at the point $B$, by the concavity
of $g$ one gets
$$\ds\f{\sin x}{x}\le \ds\f{4}{\pi }-\ds\f{4}{\pi ^2}\cdot x,
\eqno(5)$$
which is a counterpart of (4) (see \cite{54-13} for details).
Such inequalities will be used in the next sections.

\section{A counterpart of Jordan's inequality}
Another inequality, which may be found also in \cite{54-13} (see relation (18)) is the following:
$$\ds\f{\tan\ds\f{x}{2}}{x}\le \ds\f{2}{\pi },\q 0\le x\le \ds\f{\pi }{2}.
\eqno(6)$$

Here the left side is interpreted in $x=0$ as
$\lim\limits_{x\to 0}\ds\f{\tan\ds\f{x}{2}}{x}=\ds\f{1}{2}$.

For a new proof of relation (6), let us remark that, the graph of convex function
$x\to \tan\ds\f{x}{2}$ $\left(0\le x\le \ds\f{\pi }{2}\right)$ is below the segment line passing through
the points $(0,0)$ and $\left(\ds\f{\pi }{2},1\right)$.

In what follows, (6) will be called as "a counterpart of Jordan's inequality".

Our first result shows that this inequality refines Kober's inequality:

{\bf Theorem 1.}
{\it One has
$$\ds\f{1-\cos x}{x}\le \ds\f{\tan\ds\f{x}{2}}{x}\le \ds\f{2}{\pi },\q
0<x\le \ds\f{\pi }{2}.
\eqno(7)$$

}

{\bf Proof.}
By
$1-\cos x=2\sin^2 \ds\f{x}{2}$ and $2\sin\ds\f{x}{2}\cos\ds\f{x}{2}=\sin x$,
the inequality
$1-\cos x\le\tan \ds\f{x}{2}$
becomes
$\sin x\le 1$, which is true, with equality only for $x=\ds\f{\pi }{2}$.

Another proof can be obtained by letting $\tan\ds\f{x}{2}=t$, and using the formula
$\cos x=\ds\f{1-t^2}{1+t^2}$.
Then we have to prove that,
$$\ds\f{2t^2}{1+t^2}\le t\q\mbox{or}\q 2t\le 1+t^2,$$
which is $(t-1)^2\ge 0$, etc.

Now we will obtain a lower bound for $\left(\tan \ds\f{x}{2}\right)/x$:

{\bf Theorem 2.}
$$\ds\f{1}{\pi -x}\le \ds\f{\tan\ds\f{x}{2}}{x}\le \ds\f{2}{\pi },\q
0\le x\le \ds\f{\pi }{2}.
\eqno(8)$$

{\bf Proof.}
We will obtain a method of proof of (3), suggested in Section~1.
Put
$\ds\f{\pi }{2}-x$ in place of $x$ in inequality (6).
As
$$\tan\left(\left(\ds\f{\pi }{2}-\ds\f{x}{2}\right)/2\right)
=\tan\left(\ds\f{\pi }{4}-\ds\f{x}{2}\right)=\ds\f{1-\tan\ds\f{x}{2}}{1+\tan\ds\f{x}{2}},$$
after some elementary transformations we get the left side of (8).

{\bf Remark 1.}
It can be proved immediately that, the lower bounds of (7), resp. (8) cannot be compared,
i.e. some of the inequalities
$$\ds\f{1-\cos x}{x}<\ds\f{1}{\pi -x}
\q\mbox{and}\q
\ds\f{1-\cos x}{x}>\ds\f{1}{\pi -x}$$
is true for all $x\in \left(0,\ds\f{\pi }{2}\right)$.

Strong improvements of the left side of (8) will be obtained by another methods
(see Theorem 3).

Though $\ds\f{1-\cos x}{x}$ and $\ds\f{1-\cos x}{x^2}$
cannot be compared for all $x\in \left(0,\ds\f{\pi }{2}\right)$
(only for $0<x\le 1$ or $1\le x<\ds\f{\pi }{2}$), the later one is also a lower bound for
$(\tan x/2)/x$.
More generally, the following inequalities are true.

{\bf Theorem 3.}
$$\ds\f{1}{\pi }\le \ds\f{\sin x}{2x}\le \ds\f{1-\cos x}{x^2}\le \ds\f{\sin\ds\f{x}{2}}{x}
\le \ds\f{1}{2}\le \ds\f{\tan\ds\f{x}{2}}{x}\le \ds\f{2}{\pi }
\eqno(9)$$

{\bf Proof.}
The first inequality of (9) is exactly Jordan's inequality (1).
Applying
$$1-\cos x=2\sin^2\ds\f{x}{2}
\q\mbox{and}\q
\sin x=2\sin\ds\f{x}{2}\cos\ds\f{x}{2},$$
the second inequality follows by
$\ds\f{x}{2}\le \tan\ds\f{x}{2}$,
while the third one, by
$$\sin\ds\f{x}{2}\le \ds\f{x}{2}.$$

The lower bound $\ds\f{1}{\pi }$ however, is not the best one for
$$q(x)=\ds\f{1-\cos x}{x^2},\q 0<x<\ds\f{\pi }{2}.$$

{\bf Lemma 1.}
{\it The function $q$ defined above is strictly decreasing and strictly concave.
}

{\bf Proof.}
After some elementary computations (which we omit here) one obtains
$$q'(x)=\ds\f{x\sin x+2\cos x-2}{x^3},$$
$$x^4 q''(x)=x^2\cos x-4x\sin x-6\cos x+6=p(x),$$
$$p'(x)=2\sin x-2x\cos x-x^2\sin x,\q
p''(x)=-x^2\cos x<0.$$

As $p''(x)<0$, one gets
$p'(x)<p'(0)=0$ for $x>0$, so
$p(x)<p(0)=0$,
which shows that $q''(x)<0$, i.e. $q$ is strictly concave.

By letting $r(x)=x\sin x+2\cos x-2$,
remark that $r(x)<0$ can be written also as
$$\ds\f{\sin x}{2x}\le \ds\f{1-\cos x}{x^2},$$
which is the second inequality of (9).
Thus $q$ is strictly decreasing.

{\bf Theorem 4.}
$$\ds\f{4}{\pi ^2}\le \ds\f{1-\cos x}{x^2}\le \ds\f{1}{2}\le \ds\f{1}{2}\le
\ds\f{\tan\ds\f{x}{2}}{x}\le \ds\f{2}{\pi };
\eqno(10)$$
$$\ds\f{4}{\pi ^2}\le \ds\f{1-\cos x}{x^2}\le \ds\f{4}{\pi ^2}+\ds\f{4(4-\pi )}{\pi ^3}
\left(\ds\f{\pi }{2}-x\right);
\eqno(11)$$
$$\ds\f{1}{2}-\ds\f{\pi ^2-8}{\pi ^3}x\le \ds\f{1-\cos x}{x^2}\le \ds\f{1}{2}.
\eqno(12)$$

{\bf Proof.}
The first two inequalities of (10) are consequences of
$$q\left(\ds\f{\pi }{2}\right)\le q(x)\le q(0+)=\lim\limits_{x\to 0}q(x)=\ds\f{1}{2}.$$

Now, we have essentially to prove the right side of (11), as well as the left side of (12).

We will use the method of proof of inequality (5).
By writing the tangent line to the graph of function $q$ at the point
$B\left(\ds\f{\pi }{2},\ds\f{4}{\pi ^2}\right)$, as
$q'\left(\ds\f{\pi }{2}\right)=\ds\f{4(\pi -4)}{\pi ^3}$,
the right side of (11) follows.
The line passing through the point $B$ above and the point $A\left(0,\ds\f{1}{2}\right)$
has the equation
$$y=\ds\f{1}{2}+\ds\f{8-\pi ^2}{\pi ^3}\cdot x,$$
so (12) follows, as well.

As we have seen in the Introduction, the function
$x\to \ds\f{\sin x}{x}$, $0<x<\ds\f{\pi }{2}$ is strictly decreasing and concave.
This implies at once that, similarly, the function
$p(x)=\ds\f{\sin\ds\f{x}{2}}{x}$
is strictly decreasing, and concave, too.

Since
$$p(0)=\lim\limits_{x\to 0}p(x)=\ds\f{1}{2},$$
$$p'(x)=\left(\ds\f{1}{2}x\cos \ds\f{x}{2}-\sin\ds\f{x}{2}\right)/x^2,$$
and
$p\left(\ds\f{\pi }{2}\right)=\ds\f{\sqrt 2}{\pi }$,
in a similar way, as in the proof of Theorem 4 one can deduce the following results:

{\bf Theorem 5.}
$$\ds\f{\sqrt 2}{\pi }\le \ds\f{\sin\ds\f{x}{2}}{x}\le \ds\f{\sqrt 2}{\pi }
+\ds\f{\sqrt 2(4-\pi )}{2\pi ^2}\left(\ds\f{\pi }{2}-x\right);
\eqno(13)$$
$$\ds\f{1}{2}+\ds\f{2\sqrt 2-\pi }{\pi ^2}\cdot x\le \ds\f{\sin\ds\f{x}{2}}{x}\le \ds\f{1}{2}.
\eqno(14)$$

\section[New lower and upper bounds for the counterpart of Jordan's inequality]{New lower and upper
bounds for the \\
counterpart of Jordan's inequality}

In this section we will obtain results of type (4) and (5) for the fraction $(\tan x/2)/x$.

First we need an auxiliary result:

{\bf Lemma 2.}
{\it Let $f(x)=(\tan x/2)/x$, $0<x\le \ds\f{\pi }{2}$, $f(0)=\ds\f{1}{2}$.
Then $f$ is a strictly increasing, strictly convex function.
}

{\bf Proof.}
We have to prove that $f'(x)>0$ and $f''(x)>0$ for $x\in \left(0,\ds\f{\pi }{2}\right)$.
After some elementary computations, we get
$$f'(x)=(x-\sin x)/2x^2\cos^2 x/2>0,\mbox{ as }
\sin x<x.$$

For $f''(x)$ one has
$$\left(x^2\cos^2 \ds\f{x}{2}\right)f''(x)
=(1-\cos x)-\left(1-\ds\f{\sin x}{x}\right)\left(2\cos\ds\f{x}{2}-x\tan\ds\f{x}{2}\right).$$

From the left side of inequality (2) one has
$$1-\ds\f{\sin x}{x}<\ds\f{1-\cos x}{2},$$
thus by assuming
$2\cos\ds\f{x}{2}-x\tan\ds\f{x}{2}>0$,
by
$$0<2\cos\ds\f{x}{2}-x\tan\ds\f{x}{2}<2\cos \ds\f{x}{2}<2$$
we get
$$\left(x^2\cos^2 \ds\f{x}{2}\right)f''(x)>
(1-\cos x)-\left(\ds\f{1-\cos x}{2}\right)\cdot 2=0,$$
implying $f''(x)>0$.

On the other hand, for values of $x$ such that eventually
$$2\cos \ds\f{x}{2}-x\tan\ds\f{x}{2}<0,$$
as
$1-\cos x>0$, $1-\ds\f{\sin x}{x}>0$,
clearly $f''(x)>0$.
This proves the strict convexity of $f$, too.

{\bf Remark 2.}
As $f\left(\ds\f{\pi }{2}\right)=\ds\f{2}{\pi }$,
by the strict monotonicity of $f$, a new proof of (6) follows.

{\bf Theorem 6.}
{\it For all $0\le x\le \ds\f{\pi }{2}$ one has the double-inequality
$$\ds\f{2}{\pi }\le \ds\f{2(\pi -2)}{\pi ^2}\left(x-\ds\f{\pi }{2}\right)
\le \ds\f{\tan\ds\f{x}{2}}{x}\le \ds\f{4-\pi }{\pi ^2}\cdot x+\ds\f{1}{2}.
\eqno(15)$$

}

{\bf Proof.}
As $A\left(0,\ds\f{1}{2}\right)$, $B\left(\ds\f{\pi }{2},\ds\f{2}{\pi }\right)$
are points on the graph of the convex function $f$ of Lemma 2, we can write that,
the graph of $f$ lies below of the segment $AB$ (on $\left[0,\ds\f{\pi }{2}\right]$).
Since the equation of line through $AB$ is
$$y_1(x)=\ds\f{1}{2}+\ds\f{4-\pi }{\pi ^2}\cdot x,$$
by $y_1(x)\ge f(x)$, the right side of (15) follows.

The tangent line to the graph of $f$ in the point $B$ has the equation
$$y_2(x)=\ds\f{2}{\pi }+\ds\f{2(\pi -2)}{\pi ^2}\left(x-\ds\f{\pi }{2}\right),$$
and as by convexity of $f$ one has
$f(x)\ge y_2(x)$, we get the left side of inequality (15).

{\bf Remark 3.}
As
$\ds\f{4-\pi }{\pi ^2}+\ds\f{1}{2}\le \ds\f{4-\pi }{\pi ^2}\cdot \ds\f{\pi }{2}+\ds\f{1}{2}
=\ds\f{2}{\pi }$,
the right side of (15) offers an improvement of the Jordan counterpart (6).

Now we will show that the left side of (15) gives an improvement of left side of (7),
i.e. a stronger improvement of Kober's inequality will be obtained:

{\bf Theorem 7.}
{\it One has
$$\ds\f{1-\cos x}{x}\le \ds\f{2}{\pi }+\ds\f{2(\pi -2)}{\pi ^2}\left(x-\ds\f{\pi }{2}\right),\q
0<x\le \ds\f{\pi }{2}.
\eqno(16)$$

}

{\bf Proof.}
Replacing $x$ with $\ds\f{\pi }{2}-x$ in (10), the new inequality becomes
(after some elementary computations, which we omit here)
$$\ds\f{\sin x}{x}\ge \ds\f{2}{\pi }+\ds\f{\pi -2}{\pi ^2}(\pi -2x).$$

This is exactly inequality (4) of Section 1, thus relation (16) follows.

We now state a general form of Theorem 3:

{\bf Theorem 8.}
{\it Let $r\in \left(0,\ds\f{\pi }{2}\right]$.
Then for any $x\in \left(0,\ds\f{\pi }{2}\right]$ one has
$$\ds\f{\tan \ds\f{x}{2}}{x}\le \ds\f{1}{2}
+\ds\f{\left(\tan\ds\f{r}{2}\right)/r-\ds\f{1}{2}}{r},\q
0<x\le r
\eqno(17)$$
and
$$\ds\f{\tan\ds\f{x}{2}}{x}\ge \ds\f{\tan\ds\f{r}{2}}{r}
+\left(\ds\f{r-\sin r}{2r^2\cos^2 \ds\f{r}{2}}\right)(x-r)
\eqno(18)$$

}

{\bf Proof.}
Apply the same method as in the proof of Theorem 4, by letting
$B(r,f(r))=B(r,(\tan r/2)/r)$ in place of $B(\pi /2,2/\pi )$.
Then (11) and (12) will follow on base of computations done in Lemma 1.

{\bf Remark 4.}
For $r=\ds\f{\pi }{2}$ relation (17) and (18) imply the double-inequality (15).

Inequalities of a different type may be deduced from the following auxiliary result:

{\bf Lemma 3.}
{\it Put $A(x)=\ds\f{2}{x}-\ds\f{1}{\tan\ds\f{x}{2}}$, $0<x\le \ds\f{\pi }{2}$.

Then $A$ is a strictly increasing, strictly convex function.
}

{\bf Proof.}
$A'(x)=\left(x^2-4\sin^2 \ds\f{x}{2}\right)/2x^2\sin^2 \ds\f{x}{2}>0$
by $\sin \ds\f{x}{2}<\ds\f{x}{2}$.
$$A''(x)=\left(8\sin^3 \ds\f{x}{2}-x^3\cos \ds\f{x}{2}\right)/2x^3 \sin^3 \ds\f{x}{2}>0$$
by the known inequality (due to Adamovi\'c-Mitrinovi\'c, see \cite{54-3})
$$\ds\f{\sin t}{t}>\sqrt[3]{\cos t},\q 0<t<\ds\f{\pi }{2}.
\eqno(19)$$
Thus $A$ is strictly convex, too.

We state the following result:

{\bf Theorem 9.}
{\it
$$\ds\f{1}{2}\le \ds\f{\tan\ds\f{x}{2}}{x}\le \ds\f{1}{2-\ds\f{2(4-\pi )}{\pi ^2}x^2}
\le \ds\f{1}{2-x\left(\ds\f{4-\pi }{\pi }\right)}\le \ds\f{2}{\pi },
\eqno(20)$$
where $0\le x\le \ds\f{\pi }{2}$.
}

{\bf Proof.}
Since
$A\left(\ds\f{\pi }{2}\right)=\ds\f{4}{\pi }-1$ and
$A(0):=\lim\limits_{x\to 0}A(x)=0$,
and $A$ is convex on $\left[0,\ds\f{\pi }{2}\right]$,
$A(x)\le \ds\f{2(4-\pi )}{\pi ^2}\cdot x$,
so after minor transformation we obtain the second inequality of (20).
The other relations of (20) can be verified by taking into account $0\le x\le \ds\f{\pi }{2}$.

\section{Related inequalities}

We now will apply a method of proof of (2) in paper \cite{54-9}.
This will be based on the famous Hermite-Hadamard integral inequality, as well as a generalization
obtained by the author in 1982 (\cite{54-5}, see also \cite{54-6}):

{\bf Lemma 4.}
{\it Let $f:[a,b]\to \mathbb{R}$ be a continuous convex function.
Then
$$f\left(\ds\f{a+b}{2}\right)\le \ds\f{1}{b-a}\int_a^b f(t)dt
\le \ds\f{f(a)+f(b)}{2}.
\eqno(21)$$

When $f$ is strictly convex, all inequalities in (13) are strict.
When $f$ is concave (strictly concave), then the inequalities in (21) are reversed.
}

For many applications, refinements, and generalizations of this inequality, see e.g. the book \cite{54-15}.

The following generalization of left side of (21) is due to the author:

{\bf Lemma 5.}
{\it Suppose that $f:[a,b]\to \mathbb{R}$ is $2k$-times differentiable, and
$f^{(2k)}(x)\ge 0$ for all $x\in (a,b)$.
Then
$$\int_a^b f(t)dt\ge \sum_{j=0}^{k-1}\ds\f{(b-a)^{2j+1}}{2^{2j}(2j+1)!}f^{(2j)}\left(\ds\f{a+b}{2}\right).
\eqno(22)$$

When $f^{(2k)}(x)>0$, the inequality is strict.
}

Inspired by (22), in 1989 H. Alzer \cite{54-1} proved the following counterpart:

{\bf Lemma 6.}
{\it With the same conditions as in Lemma 3, one has
$$\int_a^b f(t)dt\le \ds\f{1}{2}\sum_{i=1}^{2k-1}\ds\f{(b-a)^i}{i!}
[f^{(i-1)}(a)+(-1)^{i-1}f^{(i-1)}(b)].
\eqno(23)$$

When $f^{(2k)}(x)>0$, the inequality is strict.
}

{\bf Remark 5.}
For a common generalization of (22) and (23), see \cite{54-8}.

Particularly, applying (22) and (23) for $k=2$, along with (21), we get the following:

{\bf Lemma 7.}
{\it Suppose that $f:[0,x]\to \mathbb{R}$ is a $4$-times differentiable function such that
$f''(t)<0$ and $f^{(4)}(t)>0$.
Then one has the inequalities
$$f\left(\ds\f{x}{2}\right)+\ds\f{x^2}{24}f''\left(\ds\f{x}{2}\right)
<\ds\f{1}{x}\int_0^x f(t)dt<f\left(\ds\f{x}{2}\right)
\eqno(24)$$
and
$$\ds\f{f(0)+f(x)}{2}\le \ds\f{1}{x}\int_0^x f(t)dt$$
$$\le \ds\f{f(0)+f(x)}{2}+\ds\f{x}{4}[f'(0)-f'(x)]
+\ds\f{x^2}{12}[f''(0)+f''(x)].
\eqno(25)$$

}

{\bf Proof.}
Apply Lemma 4 for the concave function $f$ on $[a,b]=[0,x]$.
Then the right side of (24) follows.
Applying Lemma 5 for $k=2$ we get the inequality
$$\ds\f{1}{b-a}\int_a^b f(t)dt\ge f\left(\ds\f{a+b}{2}\right)
+\ds\f{(b-a)^2}{24}f''\left(\ds\f{a+b}{2}\right),
\eqno(26)$$
so the left side of (24) follows for $[a,b]=[0,x]$.

In the same manner, applying Lemma 6 for $k=2$, we get
$$\ds\f{1}{b-a}\int_a^b f(x)dx\le \ds\f{f(a)+f(b)}{2}+\ds\f{b-a}{4}[f'(a)-f'(b)]$$
$$+\ds\f{(b-a)^2}{12}[f''(a)+f''(b)]
\eqno(27)$$
so (25) will be a consequence of (27), combined with Lemma 4.

We now are in a position to deduce the following trigonometric inequalities:

{\bf Theorem 10.}
{\it For all $0<x<\ds\f{\pi }{2}$ the following inequalities are true:
$$\ds\f{\sin x}{x}<\ds\f{1-\cos x}{x}<\sin\ds\f{x}{2},
\eqno(28)$$
$$\ds\f{1}{1+\cos \ds\f{x}{2}}<\ds\f{\tan\ds\f{x}{2}}{x}
<\ds\f{1}{4}\left(1+1/\cos^2 \ds\f{x}{2}\right).
\eqno(29)$$

}

{\bf Proof.}
Apply Lemma 4 for the strictly concave function
$$f(t)=\sin t\mbox{ for } [a,b]=[0,x].$$
As
$\ds\int_0^x \sin tdt=1-\cos x,$
(28) follows.

For the proof of (29), put
$$f(t)=\ds\f{1}{2\cos^2 \ds\f{t}{2}}.$$
Since
$$2f'(t)=\sin \ds\f{t}{2}\cdot \cos^{-3} \ds\f{t}{2},$$
$$2f''(t)=\ds\f{1}{2}\cos^{-2}\ds\f{t}{2}+\ds\f{3}{2}\sin^2 \ds\f{t}{2}\cos^{-4}\ds\f{t}{2}>0,$$
$f$ will be a strictly convex function.
As
$$\int_0^x \ds\f{1}{2\cos^2 \ds\f{t}{2}}=\tan\ds\f{x}{2},$$
by (21) the double inequality (29) follows, by remarking on the left side that
$$2\cos^2 \ds\f{x}{4}=1+\cos \ds\f{x}{2}.$$

{\bf Remark 6.}
As
$\sin\ds\f{x}{2}\le \sin \ds\f{\pi }{4}=\ds\f{\sqrt 2}{2}\approx 0.7$
and
$\ds\f{2}{\pi }\approx 0.63$ the right side of (28) and (7) cannot be compared
(for all values of $x$ in $\left[0,\ds\f{\pi }{2}\right)$).
Similarly, the right side of (29) takes the greatest value
$\ds\f{3}{4}=0.75$ and the least value
$\ds\f{1}{2}=0.5$, so (29) improves (6) for certain values of $x$, and vice-versa,
(6) is strong than the right side of (29) for other values of $x$.

{\bf Theorem 11.}
{\it For all $0<x<\ds\f{\pi }{2}$ the following hold true:
$$\left(\sin \ds\f{x}{2}\right)\left(1-\ds\f{x^2}{24}\right)
<\ds\f{1-\cos x}{x}<\sin\ds\f{x}{2},
\eqno(30)$$
$$\ds\f{\sin x}{2}<\ds\f{1-\cos x}{x}<\ds\f{\sin x}{2}+\ds\f{(1-\cos x)x}{4}-\ds\f{x^2\sin x}{12}.
\eqno(31)$$

}

{\bf Proof.}
Apply (25) for the same function $f(t)=\sin t$ on $[0,x]$.
We omit the details.

{\bf Remark 7.}
Applying Lemma 6 for the function
$$f(t)=\ds\f{1}{2\cos^2 \ds\f{t}{2}},$$
further improvements of type (29) can be deduced.


\newpage
\setcounter{section}{0}

\bc
{\Large\bf 2. A Note on Certain Inequalities for Hyperbolic and Trigonometric Functions}
\ec

\begin{abstract}
The aim of this paper is to prove some improvements on certain inequalities for hyperbolic
and trigonometric functions which have been recently considered in a paper by
R. Kl\'en, M. Visuric and M. Vuorinen \cite{52-1}.
Related inequalities will be pointed out, too.
\end{abstract}

\section{Introduction}

During the last years there has been a great interest in trigonometric and hyperbolic inequalities.
For many papers on the famous Jordan's or Cusa-Huygens' Wilker's or Huygens' inequalities see
e.g. \cite{52-3}, \cite{52-4}, \cite{52-2} and the references therein.

In the recent interesting paper \cite{52-1}, the authors have considered some hyperbolic,
trigonometric or hyperbolic-trigonometric inequalities.
One of the main ingredients of the proofs in \cite{52-1} was the series representations
of trigonometric functions.
For this reason, the approximation by this method yield results e.g. on interval $(0,1)$
in  place of the natural $(0,\pi /2)$.
On the other hand, there exist results in another part of mathematics, namely the theory
of means, which applied here can lead to improvements of the obtained theorems.

We will present in this Introduction two examples from \cite{52-1}, where stronger result
can be immediately deduced.

Theorem 4.3(iii) of \cite{52-1} states that
$$\ds\f{\tanh x}{x}\le \ds\f{\tanh (kx)}{kx}
\eqno(1.1)$$
holds true for all $k\in (0,1)$ and $x\in (0,1)$.
The proof is based on the series expansion of $\tanh x$.
But we can show that, (1.1) holds for $k\in (0,1)$ and any $x>0$.

Indeed, the application
$f_1(x)=\ds\f{\tanh x}{x}$
has a derivative
$$f'_1(x)=\ds\f{x-\sinh x\cdot \cosh x}{x^2\cosh^2 x}<0$$
as $\sinh x>x$ and $\cosh x>1$ for $x>0$.
Thus the function $f_1$ is strictly decreasing, implying $f_1(x)<f_1(kx)$,
as $x>kx$ for any $x>0$, $k\in (0,1)$.
This implies (1.1), with strict inequality.

Another example will be the right side inequality of Theorem 3.1 of \cite{52-1}, namely:
$$\ds\f{\sin x}{x}\le 1-\ds\f{2x^2}{3\pi ^2},\q
0<x<\ds\f{\pi }{2}.
\eqno(1.2)$$

A much stronger inequality (in fact, the best possible one of this type) is the following:
$$\ds\f{\sin x}{x}\le 1-\ds\f{1}{3\pi }x^2,\q 0<x<\ds\f{\pi }{2}.
\eqno(1.3)$$

Consider the application
$f_2(x)=ax^3+\sin x-x$,
where
$a=\ds\f{1}{3\pi }$ and $0<x<\pi /2$.
As
$f'_2(x)=3ax^2+\cos x-1$, $f''(x)=6ax-\sin  x\le  0$
by
$\ds\f{\sin x}{x}\ge 6a=\ds\f{2}{\pi }$,
i.e. Jordan's inequality.
Therefore
$f'_2(x)<f'_2(0)=0$ and
$f_2(x)<f_2(x)=0$ for $x>0$.
This gives inequality (1.3).

Clearly
$1-\ds\f{1}{3\pi }x^2<1-\ds\f{2}{3\pi ^2}x^2$
as $\pi >2$;
so (1.3) is better than (1.2).
We note that the authors of \cite{52-1} had used (1.2) in the proof of inequality
$$\ds\f{\sin x}{x}\le \cos^3 \ds\f{x}{2},\q 0<x<\sqrt{27/5}.
\eqno(1.4)$$

Since the method is by series expansion, and use of (1.2), clearly (1.3) would offer an improvement
of $\sqrt{27/5}$ of (1.4), too.
We omit here the details.

We will consider in what follows, the following interesting inequalities from \cite{52-1}:
$$\ds\f{x^2}{\sinh^2 x}<\ds\f{\sin x}{x}<\ds\f{x}{\sinh x},\mbox{ for } x\in \left(0,\ds\f{\pi }{2}\right)
\eqno(1.5)$$
$$\ds\f{1}{\sqrt{\cosh x}}<\ds\f{x}{\sinh x}<\ds\f{1}{\sqrt[4]{\cosh x}}
\mbox{ for } x\in (0,1)
\eqno(1.6)$$
$$\ds\f{1}{\cosh x}<\ds\f{\sin x}{x}<\ds\f{x}{\sinh x}\mbox{ for } x\in \left(0,\ds\f{\pi }{2}\right)
\eqno(1.7)$$
$$\ds\f{2+\cos x}{3}\le \ds\f{x}{\sinh x}\mbox{ for } x\in (0,1)
\eqno(1.8)$$
$$\ds\f{1}{\cosh x}\le \ds\f{1+\cos x}{2}\mbox{ for } x\in \left(0,\ds\f{\pi }{2}\right)
\eqno(1.9)$$
$$\ds\f{1}{(\cos x)^{2/3}}<\cosh x<\ds\f{1}{\cos x}\mbox{ for } x\in \left(0,\ds\f{\pi }{4}\right)
\eqno(1.10)$$
and
$$\ba{rl}
{\rm (i)} & x/\arcsin x\le \sin x/x;\\
{\rm (ii)} & x/\arcsin hx\le \sinh x/x;\\
{\rm (iii)} & x/\arctan x\le \tan x/x;\\
{\rm (iv)} & x/{\rm arctanh}\, x\le \tanh x
\ea\q
x\in (0,1)
\eqno(1.11)$$

In what follows, we shall study improvements and inequalities connected with relations (1.5)-(1.10),
as well as a general result, which particularly will imply all results (1.11), even with improved intervals
[e.g. $(0/\pi /2)$ for (iii)].

\section{Means of two arguments}

The logarithmic and identric means of two positive numbers $a$ and $b$ are defined by
$$L=L(a,b)=\ds\f{a-b}{\ln a-\ln b}\mbox{ for } a\ne b,\ L(a,a)=a
\eqno(2.1)$$
and
$$I=I(a,b)=\ds\f{1}{e}(b^b/a^a)^{1/(b-a)}\mbox{ for } a\ne b,\ I(a,a)=a,
\eqno(2.2)$$
respectively.
Put
$$\ba{l}
A=A(a,b)=\ds\f{a+b}{2},\q
G=G(a,b)=\sqrt{ab},\medskip \\
Q=Q(a,b)=\sqrt{\ds\f{a^2+b^2}{2}}
\ea
\eqno(2.3)$$
for the usual arithmetic, geometric means and root square means.

The Seiffert means $P$ and $T$ are defined by
$$P=P(a,b)=\ds\f{a-b}{2\arcsin\ds\f{a-b}{a+b}},\q a\ne b,\q P(a,a)=a
\eqno(2.4)$$
and
$$T=T(a,b)=\ds\f{a-b}{2\arctan \left(\ds\f{a-b}{a+b}\right)},\q
a\ne b,\q T(a,a)=a.
\eqno(2.5)$$

There are many inequalities in the literature for these means H.-J. Seiffert proved that
(in what follows, assume $a\ne b$);
$$L<P<I
\eqno(2.6)$$
and
$$A<T<Q.
\eqno(2.7)$$

J. S\'andor proved that
$$\sqrt[3]{A^2\cdot G}<P<\ds\f{G+2A}{3}<I
\eqno(2.8)$$
and more generally,
$$\sqrt[3]{y_n^2 x_n}<P<\ds\f{x_n+2y_n}{3},
\eqno(2.9)$$
for all $n\ge 0$, where $(x_n)$ and $(y_n)$ are two sequences defined recurrently by
$$\ba{l}
x_0=G,\q
y_0=A,\medskip \\
x_{n+1}=\ds\f{x_n+y_n}{2},\q
y_{n+1}=\sqrt{x_{n+1}y_n},\q
n\ge 0.
\ea
\eqno(2.10)$$

Clearly, the first inequalities of (2.8) are given by (2.9) for $n=0$.
When $n=1$, we get the stronger relations
$$\sqrt[3]{A\left(\ds\f{A+G}{2}\right)^2}<P<
\ds\f{1}{3}\left(\ds\f{A+G}{2}+2\sqrt{A\left(\ds\f{A+G}{2}\right)}\right),
\eqno(2.11)$$
but one can obtain in fact, infinitely many improvements, by the recurrence (2.10) and inequalities (2.9).

For these and other results for the mean $P$, see the author's paper \cite{52-5}.

For the means $L$ and $I$ many inequalities can be found in \cite{52-6}, \cite{52-7}.
For example, we quote the following:
$$G<L<I<A;
\eqno(2.12)$$
$$L<A_{1/3}
\eqno(2.13)$$
(where $A_r=A_r(a,b)=\left(\ds\f{a^r+b^r}{2}\right)^{1/r}$ $(r\ne 0)$; $A_0=G$).

Note that (2.13) is due to T.P. Lin.
Two other famous inequalities of interest here are
$$\sqrt[3]{G^2\cdot A}<L<\ds\f{2G+A}{3},
\eqno(2.14)$$
where the left side is due to E.B. Leach and M.C. Sholander,
while the right side to G. P\'olya - G. Szeg\"o and B.C. Carlson
(see \cite{52-6}, \cite{52-7} for exact references).

Improvements of type (2.9) for the mean $L$ have been obtained by the author
in 1996 (see \cite{52-8} and the references therein) as follows:
$$\sqrt[3]{b_n^2a_n}<L<\ds\f{2b_n+a_n}{3},
\eqno(2.15)$$
for all $n\ge 0$,
where the sequences $(a_n)$ and $(b_n)$ are defined by the recurrences
$$\ba{l}
a_0=a,\q
b_0=b,\q
b_1=\sqrt{ab},\medskip \\
a_{n+1}=\ds\f{a_n+b_n}{2},\q
b_{n+1}=\sqrt{a_{n+1}b_n},\q
(n\ge 1).
\ea
\eqno(2.16)$$

By (2.15), the following improvements of (2.14) are obtainable:
$$\sqrt[3]{G\left(\ds\f{A+G}{2}\right)^2}<L<
\ds\f{1}{3}\left(\ds\f{A+G}{2}+2\sqrt{G\left(\ds\f{A+G}{2}\right)}\right).
\eqno(2.17)$$

This follows for $n=1$ from (2.15),
Continuing with $n=2,\ldots $ other improvements can be deduced.

The similar results for the mean $T$ are the following (see \cite{52-9})
$$\sqrt[3]{Q^2\cdot A}<T<\ds\f{A+2Q}{3}
\eqno(2.18)$$
$$\sqrt[3]{v_n^2u_n}<T<\ds\f{u_n+2v_n}{3},\q n\ge 0,
\eqno(2.19)$$
where $(u_n)$, $(v_n)$ are defined by
$$u_0=A,\q v_0=Q,\q u_{n+1}=\ds\f{u_n+v_n}{2},\q
v_{n+1}=\sqrt{u_{n+1}v_n}\q (n\ge 0).$$

Particularly, for $n=1$ from (2.19) we get
$$\sqrt[3]{Q\left(\ds\f{Q+A}{2}\right)^2}<T<
\ds\f{1}{3}\left(\ds\f{Q+A}{2}+2\sqrt{Q\left(\ds\f{Q+A}{2}\right)}\right).
\eqno(2.20)$$

We note that E. Neuman and J. S\'andor \cite{52-10} have generalized relations (2.9), (2.15) and (2.20)
to a more general mean, called as the "Schwab-Borchardt mean".

\section[Trigonometric and hyperbolic inequalities from inequalities for means]{Trigonometric
and hyperbolic \\
inequalities from inequalities for means}

Before the applications of means to trigonometric and hyperbolic inequalities,
the following identities will be pointed out.

Let $x\ne 0$.
$$L(e^x,e^{-x})=\ds\f{\sinh x}{x},\q
I(e^x,e^{-x})=e^{t\coth t-1}
\eqno(3.1)$$
and
$$G(e^x,e^{-x})=1,\q
A(e^x,e^{-x})=\cosh x,\q
Q(e^x,e^{-x})=\cosh 2x
\eqno(3.2)$$
$$P(e^x,e^{-x})=\ds\f{\sinh x}{\arcsin (\tanh x)},\q
T(e^x,e^{-x})=\ds\f{\sinh x}{\arctan (\tanh x)}
\eqno(3.3)$$
$$A_{1/3}(e^x,e^{-x})=\left(\ds\f{e^{x/3}+e^{-x/3}}{2}\right)^3
=\ds\f{\cosh x+3\cosh x/3}{4}
\eqno(3.4)$$
$$P(1+\sin x,1-\sin x)=\ds\f{\sin x}{x},\q
x\in \left(0,\ds\f{\pi }{2}\right)
\eqno(3.5)$$
$$\ba{l}
G(1+\sin x,1-\sin x)=\cos x,\q
A(1+\sin x,1-\sin x)=1,\medskip \\
Q(1+\sin x,1-\sin x)=\sqrt{1+\sin^2 x}
\ea
\eqno(3.6)$$
$$A_{1/3}(1+\sin x,1-\sin x)=\ds\f{1+3\cos x\sqrt{(1+\cos x)/2}}{4}.
\eqno(3.7)$$

For the proof of (3.4) and (3.7) use
$(a+b)^3=a^3+3ab(a+b)+b^3$
and the definitions of $\cosh x$, while for (3.7) remark that
$$\left(\sqrt[3]{1+\sin x}+\sqrt[3]{1-\sin x}\right)^3=2+3\cos x\left(\sqrt{1+\sin x}+\sqrt{1-\sin x}\right).$$

As
$$\left(\sqrt{1+\sin x}+\sqrt{1-\sin x}\right)^2=2+2\cos x,$$
(3.7) also follows.

There can be written many inequalities for $\ds\f{\sin x}{x}$ and $\ds\f{\sinh x}{x}$,
but first we state the most interesting and much studied ones.
By (2.14) and (3.1), (3.2) we get
$$\sqrt[3]{\cosh x}<\ds\f{\sinh x}{x}<\ds\f{\cosh x+2}{3},\q x\ne 0.
\eqno(3.8)$$

Here the right side of (3.8) is the "hyperbolic Cusa-Huygens" inequality (see \cite{52-2}),
while the left side was discovered (by other methods) by I. Lazarevi\'c (see \cite{52-11}).

Our inequality (2.17) improves (3.8) as follows:
$$\sqrt[3]{\left(\ds\f{\cosh x+1}{2}\right)^2}<\ds\f{\sinh x}{2}
<\ds\f{1}{3}\left(\ds\f{\cosh x+1}{2}+2\sqrt{\ds\f{\cosh x+1}{2}}\right).
\eqno(3.9)$$

Applying (2.8) with (3.2)-(3.3), we get
$$\sqrt[3]{\cosh^2 x}<\ds\f{\sinh x}{\arcsin (\tanh x)}
<\ds\f{1+2\cosh x}{3}
\eqno(3.10)$$

Similarly, from (2.18) and (3.2)-(3.3),
$$\sqrt[3]{(\cosh x)(\cosh 2x)}<\ds\f{\sinh x}{\arctan (\tanh x)}
<\ds\f{\cosh x+2\sqrt{\cosh 2x}}{3}
\eqno(3.11)$$

Applying the left side of (2.6) and (3.1), (3.3), we get
$\arcsin (\tanh x)<x$ for $x>0$.
This implies that
$$\tanh x<\sin x\mbox{ for } x\in \left(0,\ds\f{\pi }{2}\right).
\eqno(3.12)$$

We note that, in \cite{52-1} relation (3.12) is proved for $x\in (0,1)$.
Applying now (2.8) with (3.5), we get the double inequality
$$\sqrt[3]{\cos x}<\ds\f{\sin x}{x}<\ds\f{\cos x+2}{3}
\eqno(3.13)$$

The right side of (3.13) is the famous "Cusa-Huygens inequality"
(see e.g. \cite{52-2}), while the left side of this inequality is due to A. Adamovi\'c
and D.S. Mitrinovi\'c (see \cite{52-11}).

Among the many other inequalities, we mention only
$$\ds\f{\sinh x}{x}<\ds\f{\cosh x+3\cosh x/3}{4},\q x\ne 0
\eqno(3.14)$$
which follows by (3.1), (3.4) and (2.13), as well as
$$\sqrt[3]{\cos^2 x}<L(1+\sin x,1-\sin x)$$
$$<\ds\f{1+3\cos x\sqrt{(1+\cos x)/2}}{4},\q
x\in (0,\pi /2)
\eqno(3.15)$$
which is a consequence of (2.13) and (3.7), as well as (2.14) and (3.6).
Particularly, the left side of (3.15) offers the following inequality of a new type:
$$\ln\left(\ds\f{1+\sin x}{1-\sin x}\right)<\ds\f{2\sin x}{\sqrt[3]{\cos x}},\q
x\in \left(0,\ds\f{\pi }{2}\right).
\eqno(3.16)$$

\section{Refined inequalities}

First we prove that (1.8) holds true for all $x\in (0,\pi /2)$,
and that this improves the hyperbolic Cusa-Huygens inequality (3.8).

{\bf Theorem 4.1.}
$$\ds\f{\sinh x}{x}<\ds\f{3}{2+\cos x}<\ds\f{2+\cosh x}{3}
\mbox{ for } 0<x<\ds\f{\pi }{2}
\eqno(4.1)$$

{\bf Proof.}
The first inequality of (4.1) may be written as
$$f_3(x)=3x-2\sinh x-\sinh x\cdot \cos x\ge 0\mbox{ on } \left[0,\ds\f{\pi }{2}\right).$$
By a simple computation we obtain
$f''_3(x)=2(\cosh x\cdot \sin x-\sinh x)>0$
by (3.12).
This implies
$f'_3(x)\ge f'_3(0)=0$.
Thus
$f_3(x)>f_3(0)=0$ for $0<x<\pi /2$.

For the second inequality of (4.1), remark that it is equivalent with
$$f_4(x)=2\cos x+2\cosh x+\cos x\cdot \cosh x-5>0.$$

By successive differentiation we get
$$\ds\f{1}{2}f_4^{(4)}(x)=\cos x+\cosh x-2\cos x\cdot \cosh x.
\eqno(4.2)$$

Put
$f_5(x)=\cos x\cdot \cosh x-1$.
As
$$f'_5(x)=-\sin x\cdot \cosh x+\cos x\cdot \sinh x,$$
$$f''_5(x)=-2\sin x\cdot \sinh x<0
\q\mbox{for all}\q x\in \left(0,\ds\f{\pi }{2}\right),$$
one gets
$f'_5(x)<f'_5(0)=0$, so $f_5(x)\le 0$, i.e.
$$\cos x\cdot \cosh x<1\mbox{ for any } x\in \left(0,\ds\f{\pi }{2}\right).
\eqno(4.3)$$

We note that this improves the right side of (1.10).

Now, by (4.3) we get that
$$\ds\f{1}{2}f_4^{(4)}(x)>\cos x+\cosh x-2=f_6(x).$$

Since
$f'_6(x)=-\sin x+\sinh x>0$ by
$\sinh x>x>\sin x$, we get
$f_6(x)>g_6(0)=0$.
Finally, $f_4^{(4)}(x)>0$ for $x>0$, implying
$f_4^{(3)}(x)>0$, $f_4^{(2)}(x)>0$, $f'_4(x)>0$, so $f_4(x)>0$
which finishes the proof of Theorem 4.1.

As by (4.1) and (3.13) we can write
$$\ds\f{\sinh x}{x}<\ds\f{3}{2+\cos x}<\ds\f{x}{\sin x}
\eqno(4.1')$$
we have obtained a refinement of right side of (1.5).

{\bf Theorem 4.2.}
{\it
$$\left(\ds\f{x}{\sinh x}\right)^3<\ds\f{1}{\cosh x}<\ds\f{\tanh x}{x}
<\ds\f{\sin x}{x}<\ds\f{x}{\sinh x}
\eqno(4.4)$$
for all $x\in \left(0,\ds\f{\pi }{2}\right)$.
}

{\bf Proof.}
The first inequality of (4.4) follows by the left side of (3.8).
The second inequality follows by
$\ds\f{\sinh x}{x}>1$,
while the third by
$\sin x>\tanh x$, which is inequality (3.12).
The last inequality is the right side of (1.7).

{\bf Theorem 4.3.}
{\it
$$\ds\f{1}{\sqrt{\cosh x}}<\ds\f{x}{\sinh x}<\ds\f{1}{\sqrt[3]{\left(\ds\f{\cosh x+1}{2}\right)^2}}
<\ds\f{1}{\sqrt[4]{\cosh x}}
\eqno(4.5)$$
for all $x\in \left(0,\ds\f{\pi}{2}\right)$.
}

{\bf Proof.}
The right side of (4.5) is the left side of (3.9).
Now, to prove the left side inequality, consider the application
$$f_7(x)=\ds\f{\sinh x}{\sqrt{\cosh x}}-x.$$

Since
$$f'_7(x)=\ds\f{2\cosh x-\sinh^2 x-2\sqrt{\cosh x}}{2\sqrt{\cosh x}},$$
and remarking that
$$\sinh^2 x-2\cosh x+2\sqrt{\cosh x}
=\cosh^2 x-1-2\cosh x+2\sqrt{\cosh x}$$
$$=(\cosh x-1)^2+2(\sqrt{\cosh x}-1)>0$$
as $\cosh x>1$ we get that
$f'_7(x)<0$, implying
$f_7(x)<f_7(0)=0$.

Thus the left side of inequality (1.6) holds in fact for all
$x\in \left(0,\ds\f{\pi }{2}\right)$.

The last inequality for $t=\cosh x$ can be written as
$\left(\ds\f{t+1}{2}\right)^8>t^3$.
Since
$\ds\f{t+1}{2}>\sqrt t$,
this follows by $t^4>t^3$, i.e. $t>1$ which is true for $x>0$.

The following result shows  that the hyperbolic analogue of (1.4) is always true:

{\bf Theorem 4.4.}
$$\ds\f{\sinh x}{x}<\cosh^3 \ds\f{x}{3}<\ds\f{\cosh x+2}{3},\q
x>0.
\eqno(4.6)$$

{\bf Proof.}
As
$$\cosh^3 \ds\f{x}{3}=\ds\f{\cosh x+3\cosh x/3}{4},$$
the left side of (4.6) holds true by (3.14).
For the second inequality we have to prove that
$$f_8(x)=\cosh 3x-9\cosh x+8>0\mbox{ for } x>0.$$

By
$f'_8(x)=3\sinh 3x-9\sinh x$,
$f''_8(x)=9(\cosh x-\cosh x)>0$
for $x>0$ it follows
$f'_8(x)>0$, thus
$f_8(x)>0$ for $x>0$.
This finishes the proof of (4.6).

{\bf Theorem 4.5.}
{\it
$$\ds\f{8\sin x}{1+3\cos x\sqrt{(1+\cos x)/2}}
<h\left(\ds\f{1+\sin x}{1-\sin x}\right)
<\ds\f{2\sin x}{\sqrt[3]{\cos x}}
\eqno(4.7)$$
for all $x\in \left(0,\ds\f{\pi }{2}\right)$.
}

{\bf Proof.}
Apply (2.13) combined with (3.7) in order to deduce the left side of (4.2).
The right side of (2.7) coincides with (3.16).

{\bf Theorem 4.6.}
{\it
$$\ds\f{1}{\cosh x}<\ds\f{1+\cos x}{2}<\ds\f{1}{\sqrt{\cosh x}}<\ds\f{x}{\sinh x}
\q\mbox{for}\q x\in \left(0,\ds\f{\pi }{2}\right).
\eqno(4.8)$$

}

{\bf Proof.}
The first inequality of (4.3) is relation (1.9).
For the second inequality put
$$f_9(x)=(\cos x+1)\sqrt{\cosh x}-2.$$

As
$$f'_9(x)=\ds\f{f_{10}(x)}{2\sqrt{\cosh x}},$$
where
$f_{10}(x)=-2\sin x\cosh x+\cos x\sinh x+\sinh x.$

Since
$\tanh x\left(\ds\f{1+\cos x}{2}\right)<\tanh x<\sin x$
by (3.12) we get
$$(\sinh x)(1+\cos x)<2\sin x\cosh x,$$
thus
$f_{10}(x)<0$.
By $f'_9(x)<0$ one has
$f_9(x)<f_9(0)=0$.
The last inequality is the first relation of (4.5).
For the left side of (1.10), the following best possible result is true.

{\bf Theorem 4.7.}
{\it There exist a unique $\lambda \in \left(\ds\f{\pi }{4},\ds\f{\pi }{4}\right)$
such that
$$\cosh^3 \lambda \cdot \cos^2 \lambda =1.$$

For any $x\in (0,\lambda )$ one has
$\cosh^3 x\cdot \cos^2 x>1$ and any
$x\in \left(\lambda ,\ds\f{\pi }{2}\right)$ one has
$$\cosh^3 x\cdot \cos^2 x<1.
\eqno(4.9)$$

}

{\bf Proof.}
Define
$f_{11}(x)=\cosh^3 x\cdot \cos^2 x-1$.
As
$$f'_{11}(x)=\cosh^2 x\cdot \cos x(3\sinh x\cdot \cos x-2\cosh x\cdot \sin x)$$
$$=\cosh^3 x\cdot \cos^2 x\cdot f_{12}(x),$$
where
$f_{12}(x)=3\tanh x-2\tan x$.
One has
$$f'_{12}(x)=\ds\f{3\cos^2 x-2\cosh^2 (x)}{\cos^2 x\cdot \cosh^2 x}.$$

By letting
$f_{13}(x)=3\cos^2 x-2\cosh^2 x$,
clearly
$$f'_{13}(x)=-6\cos x\sin x-4\cosh x\sinh x<0,$$
so $f_{13}$ is strictly decreasing function.
One has
$f_{13}\left(\ds\f{\pi }{4}\right)<0$ as
$\cosh^2 \ds\f{\pi }{4}>1>\ds\f{3}{4}$.
Since
$f_{13}(0)=1>0$
there exists a single $x_0\in \left(0,\ds\f{\pi }{4}\right)$ such that
$f_{13}(x_0)=0$ and $f_{13}(x)>0$ for $x\in (0,x_0)$ and
$f_{13}(x)<0$ for $x\in \left(x_0,\ds\f{\pi }{2}\right)$.
This means that $f_{12}$ is strictly increasing on $(0,x_0)$ and strictly decreasing on
$(x_0,\pi /2)$.
As $\tanh \ds\f{\pi }{4}=0,655\ldots <0,66\ldots =\ds\f{2}{3}$
we get
$f_{12}\left(\ds\f{\pi }{4}\right)<0$ and also clearly
$f_{12}\left(\ds\f{\pi }{2}\right)<0$.
Thus there exists a unique $x_1\in \left(x_0,\ds\f{\pi }{4}\right)$ such that
$f_{12}(x_1)=0$ and $x_0$ is a maximum point of $f_{12}$ on $(0,x_1)$.
Thus $f_{11}$ will be strictly increasing on $(0,x_1)$ and strictly decreasing on $(x_1,\pi /2)$.
As $f_{11}(0)=0$, $f_{11}\left(\ds\f{\pi }{2}\right)<0$ and
$f_{11}\left(\ds\f{\pi }{4}\right)>0$,
there exists a unique $\lambda \in \left(\ds\f{\pi }{4},\ds\f{\pi }{2}\right)$ such that
$f_{11}(\lambda )=0$.
Thus $f_{11}$ will be strictly positive on $(0,\lambda )$ and strictly negative on
$(\lambda ,\pi /2)$.

The next theorem shows that among others, the right side of (1.10) may be improved for all
$x\in (0,\pi /2)$.

{\bf Theorem 4.8.}
{\it For all $x\in \left(0,\ds\f{\pi }{2}\right)$ one has
$$\cot x+\coth x<\ds\f{2}{x},
\eqno(4.10)$$
$$\tan x\cdot \tanh x>x^2,
\eqno(4.11)$$
$$\tan x+\tanh x>2x
\eqno(4.12)$$
$$\cos x\cdot \cosh x<\ds\f{\sin x\cdot \sinh x}{x^2}<1.
\eqno(4.13)$$

}

{\bf Proof.}
For (4.10) we have to prove that
$$f_{14}(x)=x\cos x\cdot \cosh x+x\sin x\cdot \sinh x-2\sin x\cdot \sinh x<0.$$

After simple computations we get
$$f'_{14}(x)=-\cos x\cdot \sinh x-\sin x\cdot \cosh x+2x\cdot \cos x\cdot \cosh x$$
$$=\cos x\cdot \cosh x\cdot f_{15}(x),$$
where
$f_{15}(x)=-\tan x-\tanh x+2x$.

We have
$$f'_{15}(x)=-\ds\f{1}{\cos^2 x}-\ds\f{1}{\cosh^2 x}+2$$
$$=2\cos^2 x\cosh^2 x-\cos^2 x-\cosh^2 x
<2-\cos^2 x-\cosh^2 x,$$
by $\cos x\cdot \cosh x<1$.
Since
$1-\cos^2 x=\sin^2 x$
and
$1-\cosh^2 x=-\sinh^2 x$,
$f'_{15}(x)<\sin^2 x-\sinh^2 x<0$
by $\sin x<x<\sinh x$.
Thus
$f_{15}(x)<f_{15}(0)=0$,
which proves (4.12).
Thus
$f'_{14}(x)<0$, so
$f_{14}(x)<f_{14}(0)=0$,
and inequality (4.10) is proved.

As
$\cos tx+\coth x=\ds\f{1}{\tan x}=\ds\f{1}{\tanh x}$,
remark that (4.10) may be written also as
$H(\tanh x,\tan x)>x$,
where $H(a,b)$ denotes the harmonic mean of $a$ and $b$.
As the geometric mean
$$G(\tanh x,\tan x)>H(\tanh x,\tan x)>x,$$
we get inequality (4.11).
Finally, as by (4.11)
$$\cos x\cdot \cosh x<\ds\f{\sin x\cdot \sinh x}{x^2},$$
(4.13) follows by right side of (1.5).

\section[Inequalities connecting functions with their inverses]{Inequalities connecting functions with \\
their inverses}

In this section we will consider a general result on functions and their inverses.

{\bf Theorem 5.1.}
{\it Let $f:I\to J$ be a bijective function, where $I$, $J$ are nonvoid subsets of $(0,+\infty )$.
Suppose that the function
$g(x)=\ds\f{f(x)}{x}$, $x\in I$
is strictly increasing.

Then for any $x\in I$, $y\in J$ such that $f(x)\ge y$ one has
$$f(x)f^*(y)\ge x y,
\eqno(5.1)$$
where $f^*:J\to I$ denotes the inverse function of $f$.

Under the same conditions, if $f(x)\le y$ one has the reverse inequality
$$f(x)f^*(y)\le xy.
\eqno(5.2)$$

}

{\bf Proof.}
First remark that $f$ must be strictly increasing, too.
Indeed, if $x_1,x_2\in I$ and $x_1<x_2$, then
$\ds\f{f(x_1)}{x_1}<\ds\f{f(x_2)}{x_2}$,
so
$f(x_1)<\ds\f{x_1}{x_2}f(x_2)<f(x_2)$.
Thus $f^*$ is strictly increasing, too.
Put $t=f^*(y)$.
Since $f^*$ is strictly increasing, we can write $t\le x$, so
$\ds\f{f(t)}{t}\le \ds\f{f(x)}{x}$
so $yx\le f(x)f^*(y)$, i.e. inequality (5.1) holds true.

When $y\ge f(x)$ we can write similarly $\ds\f{f(x)}{x}\le \ds\f{f(t)}{t}$,
where $t=f^*(y)\ge x$, and (5.2) follows.

Clearly one has equality in (5.1) or (5.2) only when $y=f(x)$.
The following result will be obtained with the aid of (5.1).

{\bf Theorem 5.2.}
{\it (1) For any $x\in (0,1)$ and $y\in \left(0,\ds\f{\pi }{2}\right)$
such that $y<\arcsin x$ one has
$$\arcsin x\cdot \sin y>xy.
\eqno(5.3)$$

(2) For any $x>0$, $y>0$ such that $\sinh x>y$ one has
$$\sinh x\cdot {\rm arcsinh}\, y>xy.
\eqno(5.4)$$

(3) For any $x\in \left(0,\ds\f{\pi }{2}\right)$ and $y\in (0,\infty )$ such that
$\tan x>y$ one has
$$\tan x\cdot \arctan y>xy.
\eqno(5.5)$$

(4) For any $x\in (0,1)$ and $y\in (0,\infty )$ such that ${\rm arctanh}\, x>y$ one has
$${\rm arctanh}\, x\cdot \tanh y>xy.
\eqno(5.6)$$

}

{\bf Proof.}
(1) Let $I=(0,1)$, $J=\left(0,\ds\f{\pi }{2}\right)$ and $f(x)=\arcsin x$.
Then
$$g(x)=\ds\f{\arcsin x}{x}$$
is strictly increasing, as
$$g'(x)=\ds\f{x-\sqrt{1-x^2}\arcsin x}{\sqrt{1-x^2}}>0$$
by $\arcsin x<\ds\f{x}{\sqrt{1-x^2}}$.
Indeed, by letting $x=\sin p$, this becomes
$$p<\ds\f{\sin p}{\cos p}=\tan p,$$
which is well-known.

(2) Let $I=J=(0,\infty )$ and
$f(x)=\sinh x$.
Then
$$\left(\ds\f{\sinh x}{x}\right)'=\ds\f{\cosh x\cdot x-\sinh x}{x^2}>0$$
by $\tanh x<x$, which is known.

(3) For $I=\left(0,\ds\f{\pi }{2}\right)$, $J=(0,\infty )$ and
$f(x)=\tan x$ one has
$$\left(\ds\f{f(x)}{x}\right)'=\ds\f{x-\sin x\cos x}{x^2\cos^2 x}>0$$
as $x>\sin x$ and $1>\cos x$.

(4) $I=(0,1)$, $J=(0,\infty )$,
$f(x)={\rm arctanh}\, x$.
As
$$\left(\ds\f{{\rm arctanh}\, x}{x}\right)'=\left(\ds\f{x}{1-x^2}-{\rm arctanh}\, x\right)/x^2>0$$
by ${\rm arctanh}\, x>\ds\f{x}{1-x^2}$,
which is equivalent, by letting $x=\tanh p$ by $p<\sinh p\cdot \cosh p$.
This is true, as
$\sinh p>0$ and $\cosh p>1$ for $p>0$.

As a corollary, the above theorem gives:

{\bf Theorem 5.3.}
$$\ba{rlll}
(i)  & \ds\f{x}{\arcsin x}<\ds\f{\sin x}{x} & for  & x\in (0,1)\medskip \\
(ii) & \ds\f{x}{{\rm arcsinh}\, x}<\ds\f{\sinh x}{x} & for & x>0;\medskip \\
(iii) & \ds\f{x}{\arctan x}<\ds\f{\tan x}{x} & for & x\in \left(0,\ds\f{\pi }{2}\right);\medskip \\
(iv) & \ds\f{x}{{\rm arctanh}\, x}<\ds\f{\tanh x}{x} & for & x\in (0,1).
\ea
\eqno(5.7)$$

{\bf Proof.}
(i) by (5.3) from $\arcsin x>x$ for $x\in (0,1)$.

(ii) by (5.4) from $\sinh x>x$ for $x>0$.

(iii) by (5.5) from $\tan x>x$ for $x\in \left(0,\ds\f{\pi }{2}\right)$.

(iv) by (5.6) from ${\rm arctanh}\, x>x$ for $x\in (0,1)$.

\noindent
These exact results may be compared with relations (1.11), proved in \cite{52-1}.

The following result will be an application of (5.2):

{\bf Theorem 5.4.}
$$\ba{rll}
(i) & \ds\f{x}{\arcsin x}>\ds\f{\sin\left(\ds\f{\pi }{2}x\right)}{\ds\f{\pi }{2}x}
& for\  x\in (0,1);\medskip \\
(ii) & \ds\f{\sinh x}{x}<\ds\f{x}{a\cdot \arcsin(x/a)},
& for\  x\in (0,k),\ where\\
& &  \ k>0\ and \ a=\ds\f{k}{\sinh k};\medskip \\
(iii) & \ds\f{\tan x}{x}<\ds\f{bx}{\arctan(bx)},
& for\  x\in (0,k),\ where\\
& &  \ 0<k<\pi /2\ and \
b=\ds\f{\tan k}{k};\medskip \\
(iv) & \ds\f{x}{{\rm arctanh}\, x}>\ds\f{\tanh (x/c)}{x/c},
& for\  x\in (0,k),\ where \\\
& &  \ k\in (0,1)
\ and \ c=\ds\f{{\rm arctanh}\, k}{k}.
\ea
\eqno(5.8)$$

{\bf Proof.}
(i) Since $\ds\f{\arcsin x}{x}$ is a strictly increasing function of $x$,
$$\ds\f{\arcsin x}{x}<\ds\f{\arcsin (\pi /2)}{\pi /2},$$
so
$\ds\f{2}{\pi }\arcsin x=f(x)<x$.
Now, $f$ is bijective, having the inverse
$$f^*(x)=\sin\left(\ds\f{\pi }{2}\cdot x\right).$$
Thus, relation (5.2) of Theorem 5.1, applied to $y=x$ implies (i) of (5.8).

(ii) Put $f(x)=a\cdot \sinh x<x$, and apply the same method.

(iii) Let $f(x)=\ds\f{1}{b}\cdot \tan x$ and use the monotonicity of $\ds\f{\tan x}{x}$.

(iv) Let $f(x)=c\cdot {\rm arctanh}\, x$ and use the monotonicity of (iv) of Theorem 5.3.

{\bf Remarks.}
1) When $k=1$ and $a=\ds\f{1}{\sinh 1}$ in (ii) of (5.8), we get
$$\ds\f{\sinh x}{x}<(\sinh 1)\cdot \ds\f{x}{\arcsin(\sinh 1)a}\mbox{ for } x\in (0,1).$$

2) When $k=\pi /4$ in (iii) of (5.8), we get
$$\ds\f{\tan x}{x}<\ds\f{\left(\ds\f{4}{\pi }x\right)}{\arctan \left(\ds\f{4}{\pi }x\right)}
\mbox{ for } x\in \left(0,\ds\f{\pi }{4}\right).
\eqno(5.9)$$

3) When $k=\ds\f{1}{2}$ in (iv) of (5.8), we get
$$\ds\f{x}{{\rm arctanh}\, x}>\ds\f{\tanh(x/\ln 3)}{x/\ln 3}\mbox{ for }
x\in \left(0,\ds\f{1}{2}\right).$$


\newpage
\setcounter{section}{0}

\bc
{\Large\bf 3. The Huygens and Wilker-Type Inequalities as Inequalities for Means of Two Arguments}
\ec

\begin{abstract}
The famous Huygens and Wilker type trigonometric or hyperbolic inequalities
have been extensively studied in the mathematical literature (see e.g. \cite{55-18}, \cite{55-23}, \cite{55-9}).
The special means of two arguments, as the logarithmic, the identric, or the Seiffert's means
have also been in the focus of many researchers (see e.g. \cite{55-10}, \cite{55-20}, \cite{55-13}, \cite{55-8}).
The aim of this paper is to show that the Wilker or Huygens (or Cusa-Huygens) type trigonometric
and hyperbolic inequalities can be stated in fact as certain inequalities for means of two arguments.
\end{abstract}

\noindent
{\bf AMS Subject Classification (2010):}  26D05, 26D07, 26D99.

\noindent
{\bf Keywords and phrases:}
Inequalities, trigonometric functions, hyperbolic functions, means of two arguments.

\section{Introduction}

The famous Wilker inequality for trigonometric functions states that for any $0<x<\pi /2$
one has
$$\left(\ds\f{\sin x}{x}\right)^2+\ds\f{\tan x}{x}>2,
\eqno(1)$$
while the Huygens inequality asserts that
$$\ds\f{2\sin x}{x}+\ds\f{\tan x}{x}>3.
\eqno(2)$$

Another important inequality, called as the Cusa-Huygens inequality, says that
$$\ds\f{\sin x}{x}<\ds\f{\cos x+2}{3},
\eqno(3)$$
in the same interval $(0,\pi /2)$.

The hyperbolic versions of these inequalities are
$$\left(\ds\f{\sin hx}{x}\right)^2+\ds\f{\tanh x}{x}>2;
\eqno(4)$$
$$\ds\f{2\sinh x}{x}+\ds\f{\tanh x}{x}>3;
\eqno(5)$$
$$\ds\f{\sinh x}{x}<\ds\f{\cosh x+2}{3},
\eqno(6)$$
where in all cases, $x\ne 0$.

For history of these inequalities, for interconnections between them, generalizations, etc.,
see e.g. papers \cite{55-22}, \cite{55-23}, \cite{55-10}.

Let $a,b>0$ be two positive real numbers.
The arithmetic, geometric, logarithmic and identric means of these numbers are defined by
$$A=A(a,b)=\ds\f{a+b}{2},\q
G=G(a,b)=\sqrt{ab},
\eqno(7)$$
$$L=L(a,b)=\ds\f{a-b}{\ln a-\ln b}\ (a\ne b);\q
I=I(a,b)=\ds\f{1}{e}(b^b/a^a)^{1/(b-a)}\ (a\ne b)
\eqno(8)$$
with $L(a,a)=I(a,a)=a$.

Let
$$A_k=A_k(a,b)=\left(\ds\f{a^k+b^k}{2}\right)^{1/k}\ (k\ne 0),\q A_0=G
\eqno(9)$$
be the power mean of order $k$, and put
$$Q=Q(a,b)=\sqrt{\ds\f{a^2+b^2}{2}}=A_2(a,b).
\eqno(10)$$

The Seiffert's means $P$ and $T$ are defined by
$$P=P(a,b)=\ds\f{a-b}{2\arcsin \left(\ds\f{a-b}{a+b}\right)},
\eqno(11)$$
$$T=T(a,b)=\ds\f{a-b}{2\arctan \left(\ds\f{a-b}{a+b}\right)},
\eqno(12)$$
while a new mean $M$ of Neuman and S\'andor (see \cite{55-6}, \cite{55-8}) is
$$M=M(a,b)=\ds\f{a-b}{2{\rm arcsinh} \left(\ds\f{a-b}{a+b}\right)}.
\eqno(13)$$

For history and many properties of these means we quote e.g.
\cite{55-11}, \cite{55-13}, \cite{55-6}, \cite{55-8}.

The aim of this paper is to show that inequalities of type (1)-(6) are in fact inequalities
for the above stated means.
More generally, we will point out other trigonometric or hyperbolic inequalities,
as consequences of known inequalities for means.

\section{Main results}

{\bf Theorem 1.}
{\it For all $0<a\ne b$ one has the inequalities
$$\sqrt[3]{G^2A}<L<\ds\f{2G+A}{3};
\eqno(14)$$
$$L^2A+LG^2>2AG^2
\eqno(15)$$
and
$$2LA+LG>3AG.
\eqno(16)$$

}

As a corollary, relations (4)-(6) are true.

{\bf Proof.}
The left side of (14) is a well-known inequality, due to Leach and Sholander
(see e.g. \cite{55-3}), while the right side of (14) is another famous inequality, due to
P\'olya-Steg\"o and Carlson (see \cite{55-5}, \cite{55-2}).

For the proof of (15), apply the arithmetic mean - geometric mean inequality
$$u+v>2\sqrt{uv}$$
$(u\ne v>0)$ for $u=L^2A$ and $v=LG^2$.
By the left side of (14) we get
$$L^2A+LG^2>2\sqrt{L^3AG^2}>2\sqrt{(G^2A)(AG^2)}=2AG^2,$$
and (15) follows.

Similarly, apply the inequality
$$u+u+v>3\sqrt[3]{u^2v}$$
$(u\ne v>0)$
for $u=LA$, $v=LG$.
Again, by the left side of (14) one has
$$2LA+LG>3\sqrt[3]{(L^2A^2)(LG)}>3\sqrt[3]{(G^2A)(A^2G)}=3AG,$$
and (16) follows.

Put now $a=e^x$, $b=e^{-x}$ $(x>0)$ in inequalities (14)-(16).
As in this case one has
$$A=A(e^x,e^{-x})=\ds\f{e^x+e^{-x}}{2}=\cosh x$$
$$G=G(e^x,e^{-x})=1$$
$$L=L(e^x,e^{-x})=\ds\f{e^x-e^{-x}}{2x}=\ds\f{\sinh x}{x}$$
which are consequences of the definitions (7)-(8), from (14) we get the inequalities
$$\sqrt[3]{\cosh x}<\ds\f{\sinh x}{x}<\ds\f{\cosh x+2}{3}.
\eqno(17)$$

The left side of (17) is called also as "Lazarevi\'c's inequality" (see \cite{55-5}),
while the right side is exactly the hyperbolic Cusa-Huygens inequality (6) (\cite{55-18}).

By the same method, from (15) and (16) we get the hyperbolic Wilker inequality (4),
and the hyperbolic Huygens inequality (5), respectively.

{\bf Remark 1.}
For any $a,b>0$ one can find $x>0$ and $k>0$ such that
$a=e^xk$, $b=e^{-x}k$.
Indeed, for $k=\sqrt{ab}$ and $x=\ds\f{1}{2}\ln(a/b)$ this is satisfied.
Since all the means $A,G,L$ are homogeneous of order one
(i.e. e.g. $L(kt,kp)=kL(t,p)$)
the set of inequalities (14)-(16) is in fact equivalent with the set of (4)-(6).

Thus, we could call (15) as the "Wilker inequality for the means $A,G,L$", while (16)
as the "Huygens inequality for the means $A,G,L$", etc.

The following generalizations of (15) and (16) can be proved in the same manner:

{\bf Theorem 2.}
{\it For any $t>0$ one has
$$L^{2t} A^t+L^t G^{2t}>2A^tG^{2t}
\eqno(18)$$
and
$$2L^tA^t+L^tG^t>3A^tG^t.
\eqno(19)$$

}

These will imply the following generalizations of (4) and (5) (see \cite{55-9})
$$\left(\ds\f{\sinh x}{x}\right)^{2t}+\left(\ds\f{\tanh x}{x}\right)^t>2;
\eqno(20)$$
$$2\left(\ds\f{\sinh x}{x}\right)^t+\left(\ds\f{\tanh x}{x}\right)^t>3.
\eqno(21)$$

We now state the following Cusa-Wilker, Wilker and Huygens type inequalities for the means $P,A,G$:

{\bf Theorem 3.}
{\it For all $0<a\ne b$ one has the inequalities
$$\sqrt[3]{A^2G}<P<\ds\f{2A+G}{3};
\eqno(22)$$
$$P^2G+PA^2>2GA^2,
\eqno(23)$$
and
$$2PG+PA>3AG.
\eqno(24)$$

}

As a corollary, relations (1)-(3) are true.

{\bf Proof.}
Inequalities (22) are due to S\'andor \cite{55-11}.

For the proof of (23) and (24) apply the same method as in the proof of Theorem 1, but using,
instead of the left side of (14), the left side of (22).
Now, for the proof of the second part of this theorem, put
$a=1+\sin x$, $b=1-\sin x$ $\left(x\in \left(0,\ds\f{\pi }{2}\right)\right)$
in inequalities (22)-(24).
As one has by (7) and (11)
$$A=A(1+\sin x,1-\sin x)=1$$
$$G=G(1+\sin x,1-\sin x)=\cos x$$
$$P=P(1+\sin x,1-\sin x)=\ds\f{\sin x}{x},$$
from (22) we can deduce relation (3), while (23), resp. (24) imply the classical Wilker
resp. Huygens inequalities (1), (2).

{\bf Remark 2.}
Since for any $a,b>0$ one can find $x\in \left(0,\ds\f{\pi }{2}\right)$ and $k>0$ such that
$a=(1+\sin x)k$, $b=(1-\sin x)k$
[Indeed, $k=\ds\f{a+b}{2}$, $x=\arcsin \ds\f{a-b}{a+b}$],
by the homogeneity of the means $P,A,G$ one can state that
the set of inequalities (22)-(24) is equivalent with the set (1)-(3).

Thus, e.g. inequality (23) could be called as the "classical Wilker inequality for means".

{\bf Remark 4.}
Inequalities (14) and (22) can be improved "infinitely many times" by the sequential method
discovered by S\'andor in \cite{55-16} and \cite{55-11}.
For generalization, see \cite{55-6}.

The extensions of type (18)-(19) and (20), (21) can be made here, too, but we omit further details.

We now state the corresponding inequalities for the means $T$, $A$ and $Q$
((25), along with infinitely many improvements appear in \cite{55-13}).

{\bf Theorem 4.}
{\it For all $0<a\ne b$ one has the inequalities
$$\sqrt[3]{Q^2A}<T<\ds\f{2Q+A}{3},
\eqno(25)$$
$$T^2A+TQ^2>2AQ^2
\eqno(26)$$
and
$$2TA+TQ>3AQ.
\eqno(27)$$

}

The corresponding inequalities for the means $M$, $A$ and $Q$ will be the following:

{\bf Theorem 5.}
{\it For all $0<a\ne b$ one has the following inequalities:
$$\sqrt[3]{A^2Q}<M<\ds\f{2A+Q}{3};
\eqno(28)$$
$$M^2Q+MA^2>2AQ^2;
\eqno(29)$$
$$2MQ+MA>3QA.
\eqno(30)$$

}

{\bf Proof.}
Inequalities (28) can be found essentially in \cite{55-6}.
Relations (29) and (30) can be proved in the same manner as in the preceding theorems,
by using in fact the left side of (28).

{\bf Remark 5.}
For an application, put $a=e^x$, $b=e^{-x}$ in (28).
Since
$$A=\ds\f{e^x+e^{-x}}{2}=\cosh x,\q
Q=\sqrt{\ds\f{e^{2x}+e^{-2x}}{2}}=\sqrt{\cosh (2x)}$$
and
$$a-b=2\sinh x,\q
\ds\f{a-b}{a+b}=\tanh x$$
and by
$${\rm arcsinh} (t)=\ln(t+\sqrt{t^2+1}),$$
we get
$$M(a,b)=\ds\f{\sinh x}{\ln(\tanh x+\sqrt{1+\tanh^2 x})}.$$

For example, the Wilker inequality (29) will become:
$$\ds\f{(\tanh x)^2}{\ln^2(\tanh x+\sqrt{1+\tanh^2 x})}
+\ds\f{\sinh x}{\sqrt{\cosh 2x}}\cdot \ds\f{1}{\ln(\tanh x+\sqrt{1+\tanh^2 x})}>2.
\eqno(31)$$

Since
$\sqrt{\cosh 2x}=\sqrt{\cosh^2 x-\sinh ^2 x}<\cosh x$,
here we have
$$\ds\f{\sinh x}{\sqrt{\cosh x}}>\tanh x,$$
so formula (31) is a little "stronger" than e.g. the classical form (1).

Finally, we point out e.g. certain hyperbolic inequalities, which will be the consequences of
the various existing inequalities between means of two arguments.

{\bf Theorem 6.}
$$1<\ds\f{\sinh t}{t}<e^{t\coth t-1}<\cosh t
\eqno(32)$$
$$e^{(t\coth t-1)/2}<\ds\f{\sinh t}{t}<\ds\f{\cosh t+3\cosh t/3}{4}
\eqno(33)$$
$$\sqrt[3]{\cosh t}<e^{(t\coth t-1)/2}<\ds\f{\sinh t}{t}<L(\cosh t,1)
<\left(\ds\f{\sqrt[3]{\cosh t}+1}{2}\right)^3
\eqno(34)$$
$$\sqrt[3]{\cosh t}<\ds\f{2\cosh t+1}{3}<e^{(t\coth t-1)/2}<\ds\f{\sinh t}{t}$$
$$<\ds\f{\cosh t+2}{3}<\ds\f{2\cosh t+1}{3}<e^{t\coth t-1}
\eqno(35)$$
$$\ds\f{\cosh 2t+3\cosh 2t/3}{4}<e^{2t\coth t-2}
<\ds\f{2\cosh^2 t+1}{3}
\eqno(36)$$
$$\sqrt[3]{\cosh^2 t}<P(e^t,e^{-t})
<\ds\f{2\cosh t+1}{3}
\eqno(37)$$
$$\ds\f{2}{e}\cosh t<e^{t\coth t-1}<\ds\f{2}{e}(\cosh t+1)
\eqno(38)$$
$$2\cosh^2 t-1<e^{2t\tanh t}
\eqno(39)$$
$$4\ln(\cos ht)>t \tanh t+3t\coth t-3.
\eqno(40)$$

{\bf Proof.}
For the identric mean of (8) one has
$I(e^t,e^{-t})=e^{t\coth t-1}$, so for the proof of (32) apply the known inequalities
(see the references in \cite{55-10})
$$1<L<I<A.
\eqno(41)$$
For the proof of (33) apply
$$\sqrt{G\cdot I}<L<A_{1/3}.
\eqno(42)$$
The left side of (42) is due to Alzer ([?]), while the right side to T.P. Lin (\cite{55-4}).
As
$$A_{1/3}(e^t,e^{-t})=\left(\ds\f{e^{t/3}+e^{-t/3}}{2}\right)^3
=\ds\f{e^t+e^{-t}+3(e^{2t/3}+e^{-2t/3})}{8},$$
(33) follows.
For the proof of (34) apply
$$\sqrt[3]{A\cdot G^2}<\sqrt{I\cdot G}<L<L(A,G)
\eqno(43)$$
and
$L(t,1)<A_{1/3}(t,1)$.

The first two inequalities of (43) are due to S\'andor \cite{55-15} and Alzer \cite{55-1}
respectively, while the last one to Neuman and S\'andor \cite{55-7}.
Inequality (35) follows by
$$\sqrt[3]{A^2\cdot G}<\sqrt{I\cdot G}<L<\ds\f{A+2G}{3}<I
\eqno(44)$$
and these can be found in \cite{55-11}.

For inequality (3) apply
$$A_{2/3}^2<I^2<\ds\f{2A^2+G^2}{3}.
\eqno(45)$$

The left side of (45) is due to Stolarsky \cite{55-21}, while the right side of
S\'andor-Trif \cite{55-12}.

Inequality (37) follows by (22), while (38) by
$$\ds\f{2}{e}A<I<\ds\f{2}{e}(A+G),
\eqno(46)$$
see Neuman-S\'andor \cite{55-7}.

Finally, for inequalities (39) and (40) we will use the mean $S$ defined by
$S(a,b)=(a^a\cdot b^b)^{1/(a+b)}$,
and remarking that
$S(e^t,e^{-t})=e^{\tanh t}$,
apply the following inequalities:
$$2A^2-G^2<S^2
\eqno(47)$$
(see S\'andor-Ra\c{s}a \cite{55-17}),
while for the proof of (40) apply the inequality
$$S<A^4/I^3
\eqno(48)$$
due to S\'andor \cite{55-14}.


\newpage
\setcounter{section}{0}

\bc
{\Large\bf 4. On Cusa-Huygens Type Trigonometric and Hyperbolic Inequalities}
\ec

\begin{abstract}
Recently a trigonometric inequality by N. Cusa and C. Huygens (see e.g. \cite{56-1}, \cite{56-6})
has been discussed extensively in mathematical literature (see e.g. \cite{56-6}, \cite{56-7}, \cite{56-4}).
By using a unified method based on monotonicity or convexity of certain functions,
we shall obtain new Cusa-Huygens type inequalities.
Hyperbolic versions will be pointed out, too.
\end{abstract}

\noindent
{\bf AMS Subject Classification (2010):}  26D05, 26D07.

\noindent
{\bf Keywords and phrases:}
Inequalities, trigonometric functions, hyperbolic functions.

\section{Introduction}

In recent years the trigonometric inequality
$$\ds\f{\sin x}{x}<\ds\f{\cos x+2}{3},\q 0<x<\ds\f{\pi }{2}
\eqno(1)$$
among with other inequalities, has attracted attention of several researchers.
This inequality is due to N. Cusa and C. Huygens (see \cite{56-6} for more details regarding
this result).

Recently, E. Neuman and J. S\'andor \cite{56-4} have shown that inequality (1) implies
a result by S. Wu and H. Srivastava \cite{56-10}, namely
$$\left(\ds\f{x}{\sin x}\right)^2+\ds\f{x}{\tan x}>2,\q 0<x<\ds\f{\pi }{2}
\eqno(2)$$
called as "the second Wilker inequality".
Relation (2) implies in turn the classical and famous Wilker inequality (see \cite{56-9}):
$$\left(\ds\f{\sin x}{x}\right)^2+\ds\f{\tan x}{x}>2.
\eqno(3)$$

For many papers, and refinements of (2) and (3), see \cite{56-4} and the references therein.

A hyperbolic counterpart of (1) has been obtained in \cite{56-4}:
$$\ds\f{\sinh x}{x}<\ds\f{\cosh x+2}{3},\q x>0.
\eqno(4)$$

We will call (4) as the hyperbolic Cusa-Huygens inequality, and remark that if (4)
is true, then holds clearly also for $x<0$.

In what follows, we will obtain new proofs of (1) and (4), as well as new inequalities
or counterparts of these relations.

\section{Main results}

{\bf Theorem 1.}
{\it Let $f(x)=\ds\f{x(2+\cos x)}{\sin x}$, $0<x<\ds\f{\pi }{2}$.
Then $f$ is a strictly increasing function.
Particularly, one has
$$\ds\f{2+\cos x}{\pi }<\ds\f{\sin x}{x}<\ds\f{2+\cos x}{3},\q
0<x<\ds\f{\pi }{2}.
\eqno(5)$$

}

{\bf Theorem 2.}
{\it Let $g(x)=\ds\f{x\left(\ds\f{4}{\pi } +\cos x\right)}{\sin x}$, $0<x<\ds\f{\pi }{2}$.
Then $g$ is a strictly decreasing function.
Particularly, one has
$$\ds\f{1+\cos x}{2}<\ds\f{\ds\f{4}{\pi }+\cos x}{\ds\f{4}{\pi }+1}
<\ds\f{\sin x}{x}<\ds\f{\ds\f{4}{\pi }+\cos x}{2}.
\eqno(6)$$

}

{\bf Proof.}
We shall give a common proofs of Theorems 1 and 2.
Let us define the application
$$f_a(x)=\ds\f{x(a+\cos x)}{\sin x},\q 0<x<\ds\f{\pi }{2}.$$

Then, easy computations yield that
$$\sin^2 x\cdot f'_a(x)=a \sin x+\sin x \cos x-ax\cos x-x=h(x).
\eqno(7)$$

The function $h$ is defined on $\left[0,\ds\f{\pi }{2}\right]$.
We get
$$h'(x)=(\sin x)(ax-2\sin x).$$

Therefore, one obtains that

(i) If $\ds\f{\sin x}{x}<\ds\f{a}{2}$, then $h'(x)>0$.
Thus by (7) one has $h(x)>h(0)=0$, implying $f'_a(x)>0$, i.e. $f_a$ is strictly increasing.

(ii) If $\ds\f{\sin x}{x}>\ds\f{a}{2}$, then $h'(x)<0$, implying as above that $f_a$
is strictly decreasing.

Select now $a=2$ in (i).
Then $f_a(x)=f(x)$, and the function $f$ in Theorem 1 will be strictly increasing.
Selecting $a=\ds\f{4}{\pi }$ in (ii), by the famous Jordan inequality
(see e.g. \cite{56-3}, \cite{56-7}, \cite{56-8}, \cite{56-2})
$$\ds\f{\sin x}{x}>\ds\f{2}{\pi },
\eqno(8)$$
so $f_a(x)=g(x)$ of Theorem 2 will be strictly decreasing.

Now remarking that $f(0)<f(x)<f\left(\ds\f{\pi }{2}\right)$ and
$g(0)>g(x)>g\left(\ds\f{\pi }{2}\right)$,
after some elementary transformations, we obtain relations (5) and (6).

{\bf Remarks.}
1. The right side of (5) is the Cusa-Huygens inequality (1), while the left side
seems to be new.

2. The first inequality of (6) follows by an easy computation, based on
$0<\cos x<1$.
The inequality
$$\ds\f{1+\cos x}{2}<\ds\f{\sin x}{x}
\eqno(9)$$
appeared in our paper \cite{56-5}, and rediscovered by other authors
(see e.g. \cite{56-2}).

3. It is easy to see that inequalities (5) and (6) are not comparable, i.e. none
of these inequalities implies the other one for all $0<x<\pi /2$.

Before turning to the hyperbolic case, the following auxiliary result will be proved:

{\bf Lemma 1.}
{\it For all $x\ge 0$ one has the inequalities
$$\cos x\cosh x\le 1
\eqno(10)$$
and
$$\sin x\sinh x\le x^2.
\eqno(11)$$

}

{\bf Proof.}
Let $m(x)=\cos x\cosh x-1$, $x\ge 0$.
Then
$$m'(x)=-\sin x\cosh x+\cosh x\sinh x,$$
$$m''(x)=-2\sin x\sinh x<0.$$
Thus
$m'(x)<m'(0)=0$ and $m(x)<m(0)=0$ for $x>0$, implying (10), with equality only for $x=0$.

For the proof of (11), let
$$n(x)=x^2-\sin x\sinh x.$$
Then
$$n'(x)=2x-\cos x\sinh x-\sin x\cosh x,$$
$$n''(x)=2(1-\cos x\cosh x)<0$$
by (10), for $x>0$.
This easily implies (11).

{\bf Theorem 3.}
{\it Let
$F(x)=\ds\f{x(2+\cosh x)}{\sinh x}$, $x>0$.
Then $F$ is a strictly increasing function.
Particularly, one has inequality (4).
On the other hand,
$$\ds\f{2+\cosh x}{k^*}<\ds\f{\sinh x}{x}<\ds\f{2+\cosh x}{3},\q
0<x<\ds\f{\pi }{2}
\eqno(12)$$
where
$k^*=\ds\f{\pi }{2}(2+\cosh \pi /2)/\sinh(\pi /2)$.
}

{\bf Theorem 4.}
{\it Let
$G(x)=\ds\f{x(\pi +\cosh x)}{\sinh x}$, $x>0$.
Then $G$ is a strictly decreasing function for $0<x<\pi /2$.
Particularly, one has
$$\ds\f{\pi +\cosh x}{\pi +1}<\ds\f{\sinh x}{x}<\ds\f{\pi +\cosh x}{k},\q
0<x<\ds\f{\pi }{2}
\eqno(13)$$
where
$k=\ds\f{\pi }{2}(\pi +\cosh \pi /2)/\sinh (\pi /2)$.
}

{\bf Proof.}
We shall deduce a common proofs to Theorem 3 and 4.
Put
$$F_a(x)=\ds\f{x(a+\cosh x)}{\sinh x},\q x>0.$$
An easy computation gives
$$(\sinh x)^2 F'_a(x)=g_a(x)=a\sinh x+\cosh x\sinh x-ax\cosh x-x.$$
The function $g_a$ is defined for $x\ge 0$.
As
$$g'_a(x)=(\sinh x)(2\sinh x-ax),$$
we get that:

(i) If $\ds\f{\sinh x}{x}>\ds\f{a}{2}$, then $g'_a(x)>0$.
This in turn will imply
$F'_a(x)>0$ for $x>0$.

(ii) If $\ds\f{\sinh x}{x}<\ds\f{a}{2}$, then
$F'_a(x)<0$ for $x>0$.

By letting $a=2$, by the known inequality $\sinh x>x$, we obtain the monotonicity if
$F_2(x)=F(x)$ of Theorem 3.
Since
$F(0)=\lim\limits_{x\to 0+}F(x)=3$,
inequality (4), and the right side of (12) follows.
Now, the left side of (12) follows by
$F(x)<F(\pi /2)$ for $x<\pi /2$.

By letting $a=\pi $ in (ii) we can deduce the results of Theorem 4.
Indeed, by relation (110 of the Lemma 1 one can write
$\ds\f{\sinh x}{x}<\ds\f{x}{\sin x}$
and by Jordan's inequality (8), we get
$\ds\f{\sinh x}{x}<\ds\f{\pi }{2}$ thus $a=\pi $ may be selected.
Remarking that $g(0)>g(x)>g\left(\ds\f{\pi }{2}\right)$, inequalities (13) will follow.

{\bf Remark.}
By combining (12) and (13), we can deduce that:
$$3<k^*<k<\pi +1.
\eqno(14)$$

Now, the following convexity result will be used:

{\bf Lemma 2.}
{\it Let $k(x)=\ds\f{1}{\tanh x}-\ds\f{1}{x}$, $x>0$.
Then $k$ is a strictly increasing, concave function.
}

{\bf Proof.}
Simple computations give
$$k'(x)=\ds\f{1}{x^2}-\ds\f{1}{(\sinh x)^2}>0$$
and
$$k''(x)=\ds\f{2[x^3\cosh x-(\sinh x)^3]}{x^3(\sinh x)^3}<0,$$
since by a result of I. Lazarevi\'c (see e.g. \cite{56-3}, \cite{56-4}) one has
$$\ds\f{\sinh x}{x}>(\cosh x)^{1/3}.
\eqno(15)$$
This proves Lemma 2.

{\bf Theorem 5.}
{\it Let the function $k(x)$ be defined as in Lemma 2.
Then one has
$$\ds\f{1+x^2\cdot \ds\f{k(r)}{r}}{\cosh x}\le \ds\f{x}{\sinh x}
\mbox{ for any } 0<x\le r
\eqno(16)$$
and
$$\ds\f{x}{\sinh x}\le \ds\f{1+k(r)x+k'(r)x(x-r)}{\cosh x}\mbox{ for any } 0<x,r.
\eqno(17)$$

In both inequalities (16) and (17) there is equality only for $x=r$.
}

{\bf Proof.}
Remark that $k(0+)=\lim\limits_{x\to 0+}k(x)=0$,
and that by the concavity of $k$, the graph of function $k$ is above the line segment joining
the points $A(0,0)$ and $B(r,k(r))$.
Thus
$k(x)\ge \ds\f{k(r)}{r}\cdot x$ for any $x\in (0,r]$.
By multiplying with $x$ this inequality, after some transformations, we obtain (16).

For the proof of (17), write the tangent line to the graph of function $k$ at the point $B(r,k(r))$.
Since the equation of this line is
$y=k(r)+k'(r)(x-r)$
and writing that $y\le k(x)$ for any $x>0$, $r>0$,
after elementary transformations, we get relation (17).

For example, when $r=1$ we get:
$$\left[x^2\left(\ds\f{2}{e^2-1}\right)+1\right]/\cosh x\le \ds\f{x}{\sinh x}
\mbox{ for all } 0<x\le 1
\eqno(18)$$
and
$$\ds\f{x}{\sinh x}\le \left[1+\left(\ds\f{2}{e^2-1}\right)x+
\left(\ds\f{e^4-6e^2+1}{e^4-2e^2+1}\right)x(x-1)\right]/\cosh x
\eqno(19)$$
for any $x>0.$

In both inequalities (18) and (19) there is equality only for $x=1$.

Finally, the following convexity result, with applications will be proved:

{\bf Lemma 3.}
{\it Let $j(x)=3x-2\sinh x-\sinh x\cos x$, $0<x<\ds\f{\pi }{2}$.
Then $j$ is a strictly convex function.
}

{\bf Proof.}
Since
$j''(x)=2(\cosh x\sin x-\sinh x)>0$
is equivalent to
$$\sin x>\tanh x,\q 0<x<\ds\f{\pi }{2}
\eqno(20)$$
we will show that inequality (20) holds true for any $x\in \left(0,\ds\f{\pi }{2}\right)$.
We note that in \cite{56-2} it is shown that (20) holds for $x\in (0,1)$,
but here we shall prove with another method the stronger result (20).

Inequality (20) may be written also as
$$p(x)=(e^x+e^{-x})\sin x-(e^x-e^{-x})>0.$$

Since $p''(x)=(e^x-e^{-x})(2\cos x-1)$ and $e^x-e^{-x}>0$,
the sign of $p''(x)$ depends on the sign of $2\cos x-1$.
Let $x_0\in \left(0,\ds\f{\pi }{2}\right)$ be the unique number such that
$2\cos x_0-1=0$.
Here $x_0=\arccos \left(\ds\f{1}{2}\right)\approx 1.0471$.
Thus, $\cos x$ being a decreasing function, for all $x<x_0$ one has
$\cos x>\ds\f{1}{2}$, i.e. $p''(x)>0$ in $(0,x_0)$.
This implies $p'(x)>p'(0)=0$, where
$$p'(x)=(e^x-e^{-x})\sin x+(e^x+e^{-x})\cos x-(e^x+e^{-x}).$$

This in turn gives $p(x)>p(0)=0$.

Let now $x_0<x<\pi /2$.
Then, as $p'(x_0)>0$ and $p'\left(\ds\f{\pi }{2}\right)<0$
and $p'$ being continuous and decreasing, there exists a single
$x_0<x_1<\pi /2$ such that $p'(x_1)=0$.
Then $p'$ will be positive on $(x_0,x_1)$ and negative on $\left(x_1,\ds\f{\pi }{2}\right)$.
Thus $p$ will be strictly decreasing on $\left(x_1,\ds\f{\pi }{2}\right)$, i.e.
$p(x)>p\left(\ds\f{\pi }{2}\right)>0$.
This means that, for any $x\in \left(0,\ds\f{\pi }{2}\right)$ one has
$p(x)>0$, completing the proof of (20).

Now, via inequality (1), the following improvement of relation (11) will be proved:

{\bf Theorem 6.}
{\it For any $x\in \left(0,\ds\f{\pi }{2}\right)$ one has
$$\ds\f{\sin x}{x}<\ds\f{\cos x+2}{3}<\ds\f{x}{\sinh x}.
\eqno(21)$$

}

{\bf Proof.}
The first inequality of (21) is the Cusa-Huygens inequality (1).
The second inequality of (21) may be written as $j(x)>0$, where $j$ is the function defined
in Lemma 3.
As $j'(0)=0$ and $j'(x)$ is strictly increasing, $j'(x)>0$, implying $j(x)>j(0)=0$.
This finishes the proof of (21).


\newpage
\setcounter{section}{0}

\bc
{\Large\bf 5. Note on Huygens' Trigonometric and Hyperbolic Inequality}
\ec

\begin{abstract}
By using the theory of means, we find improvements of Huygens' trigonometric inequality
$2\sin x+\tan x>3x$ for $x\in (0,\pi /2)$,
as well as its hyperbolic version
$2\sinh x+\tanh x>3x$ for $x>0$.
\end{abstract}

\noindent
{\bf AMS Subject Classification:}
26D05, 26D07, 26D99,

\noindent
{\bf Keywords and phrases:}
Inequalities,trigonometric functions, hyperbolic
functions, means of two arguments.

\section{Introduction}
The famous Huygens trigonometric inequality (see e.g. [1], [9], [4]) states that
for all $x\in \left(0,\ds\f{\pi }{2}\right)$
one has
$$2\sin x+\tan x>3x.
\eqno(1.1)$$

The hyperbolic counterpart of (1.1) has been established recently by E. Neuman
and J. S\'andor [5]:
$$2\sinh x+\tanh x>3x,\mbox{ for } x>0.
\eqno(1.2)$$

Let $a,b > 0$ be two positive real numbers.
Denote by
$$A=A(a,b)=\ds\f{a+b}{2},\q
G=G(a,b)=\sqrt{ab},\q
H=H(a,b)=\ds\f{2ab}{a+b}$$
the classical arithmetic, geometric, resp. harmonic means of $a$ and $b$, and let
$$\ba{ll}
L=L(a,b)=\ds\f{a-b}{\ln a-\ln b}\ (a\ne b), & L(a,a)=a;\medskip \\
P=P(a,b)=\ds\f{a-b}{2\arcsin \ds\f{a-b}{a+b}}\ (a\ne b), & P(a,a)=a,
\ea
\eqno(1.3)$$
be the logarithmic, resp. Seiffert's mean of $a$ and $b$.
These means have been also in the focus of many research papers in the last decades.
For a survey of results, see [8].
The main idea of this paper is that remarking that by
$$\ba{l}
P(1+\sin x,1-\sin x)=\ds\f{\sin x}{x},\medskip \\
G(1+\sin x,1-\sin x)=\cos x,\
A(1+\sin x,1-\sin x)=1
\ea
\eqno(1.4)$$
inequality (1.1) of Huygens may be written also as
$$P>\ds\f{3AG}{2G+A},
\eqno(1.5)$$
and on the other hand, as
$$\ba{l}
L(e^x,e^{-x})=\ds\f{\sinh x}{x},\medskip \\
G(e^x,e^{-x})=1,\ A(e^x,e^{-x})=\cosh x
\ea
\eqno(1.6)$$
inequality (1.2) may be written as
$$L>\ds\f{3AG}{2A+G}.
\eqno(1.7)$$
Clearly, in (1.5)
$P=P(a,b)$ with $a=1+\sin x$, $b=1-\sin x$, $x\in \left(0,\ds\f{\pi }{2}\right)$ etc.;
while in (1.7)
$L=L(a,b)$, with $a=e^x$, $b=e^{-x}$, $x>0$, etc.;
but we will show that (1.5) and (1.7) hold true for any $a,b>0$, $a\ne b$
and in fact strong improvements hold true.
Similar ideas have been exploited in [10].

\section{Main results}
As
$3AG/(G+A)=H(G,A,A)$ and
$3AG/(2A+G)=H(G,G,A)$ where
$$H(a,b,c)=3/\left(\ds\f{1}{a}+\ds\f{1}{b}+\ds\f{1}{c}\right)$$
denotes the harmonic mean of three number, the following result may be written shortly as:

{\bf Theorem 2.1.}
{\it For any $a,b>0$, $a\ne b$ one has
$$P>H(L,A)>H(G,A,A)
\eqno(2.1)$$
and
$$L>H(P,G)>H(G,G,A).
\eqno(2.2)$$

}

{\bf Proof.}
The inequalities $P>H(L,A)$ and $L>H(P,G)$
have been proved in paper [4] (see Corollary 3.2).

Now, the second inequality of (2.1), after simpliffcations becomes in fact
relation (1.7), while, similarly the second inequality of (2.2) becomes relation (1.5).

Now, we will prove inequalities (1.5) and (1.7) in the following improved forms:
$$P>\sqrt[3]{A^2G}>\ds\f{3AG}{2G+A},
\eqno(2.3)$$
resp.
$$L>\sqrt[3]{G^2A}>\ds\f{3AG}{2A+G}.
\eqno(2.4)$$

The first inequality of (2.3) is proved in [8], while the first inequality of (2.4)
is a famous inequality due to E. B. Leach and M. C. Sholander [2] (see [6] for
many other references).

Now, both right side inequalities of (2.3) and (2.4) can be written as
$$\ds\f{2t+1}{3}>\sqrt[3]{t^2},\q t>0,$$
where $t=\ds\f{G}{A}$ in (2.3), while
$t=\ds\f{A}{G}$ in (2.4).

By the arithmetic mean - geometric mean inequality one has
$$\ds\f{2t+1}{3}=\ds\f{t+t+1}{3}>\sqrt[3]{t\cdot t\cdot 1}=\sqrt[3]{t^2},$$
so (2.5) holds true for any $t>0$, $t\ne 1$.
Since $a\ne b$, in both cases $A\ne B$;
thus the proofs of (2.3) and (2.4) is finished.
This completes the proof of Theorem 2.1.

{\bf Remark 2.1.}
Both inequalities (2.1) and (2.3) offer improvements of the classical Huygens inequality (1.1).
By using (1.4) we get form (2.3):
$$\ds\f{\sin x}{x}>\sqrt[3]{\cos x}>\ds\f{3\cos x}{2\cos x+1},\q
x\in \left(0,\ds\f{\pi }{2}\right).
\eqno(2.6)$$

By using (2.4) with (1.6), we get
$$\ds\f{\sinh x}{x}>\sqrt[3]{\cosh x}>\ds\f{3\cosh x}{2\cosh x+1},\q x>0.
\eqno(2.7)$$

The first inequality of (2.6) has been discovered also by D. D. Adamovi\'c and
D. S. Mitrinovi\'c (see [3]), while the first inequality of (2.7) by I. Lazarevi\'c (see [3]).

By the similar method, (2.1) and (2.2) become:
$$\ds\f{\sin x}{x}>\ds\f{2L^*}{L^*+1}>\ds\f{3\cos x}{2\cos x+1},\q x\in \left(0,\ds\f{\pi }{2}\right),
\eqno(2.8)$$
and
$$\ds\f{\sinh x}{x}>\ds\f{2P^*}{P^*+1}>\ds\f{3\cosh x}{2\cosh x+1},\q x>0,
\eqno(2.9)$$
where
$L^*=L(1+\sin x,1-\sin x)$, $P^*=P(e^x,e^{-x})$.

By considering inequalities (2.1) and (2.3), as well as (2.2) and (2.4), one
may ask if there is a possibility of comparison of these relations.
The answer to this question will be provided by the following:

{\bf Theorem 2.2.}
{\it For any $a,b>0$, $a\ne b$ one has
$$P>H(L,A)>\sqrt[3]{A^2G}>H(G,A,A)
\eqno(2.10)$$
and
$$L>\sqrt[3]{G^2A}>H(P,G)>H(G,G,A).
\eqno(2.11)$$

}

{\bf Proof.}
The following auxiliary result will be used:

{\bf Lemma 2.1.}
{\it For any $0<x<1$ one has
$$\ds\f{2x+1}{x+2}>\sqrt[3]{x},
\eqno(2.12)$$
and for $x>1$ one has
$$\sqrt[3]{(x+1)^2}\cdot (2\sqrt[3]{x}-1)>x\sqrt[3]{4}.
\eqno(2.13)$$

}

{\bf Proof.}
Put $x=a^3$ in (2.11).
After some transformations, the inequality becomes
$$(a+1)(a-1)^3<0,$$
which is true, as $0<a<1$.

For the proof of (2.12) remark that by
$(x+1)^2>4x$
it is sufficient to prove that
$2\sqrt[3]{x^2}-\sqrt[3]{x}>\sqrt[3]{x^2}$,
i.e.
$\sqrt[3]{x^2}>\sqrt[3]{x}$,
which is true, as $x>1$.

Now, for the proof of second inequality of (2.10) we will use the inequality
$$P<\ds\f{2A+G}{3},$$
due to the author (see [8]).
This implies
$$\ds\f{G}{P}>\ds\f{3G}{G+2A},$$
so
$$\ds\f{1}{2}\left(1+\ds\f{G}{P}\right)>\ds\f{2G+A}{G+2A}.$$

Since by (2.11)
$$\ds\f{2G+A}{G+2A}>\sqrt[3]{\ds\f{G}{A}}$$
(put $x:=\ds\f{G}{A}$), we get
$$\ds\f{1}{2}\left(1+\ds\f{G}{P}\right)>\sqrt[3]{\ds\f{G}{A}},
\mbox{ or }
\ds\f{G+P}{2GP}>\ds\f{1}{\sqrt[3]{G^2A}},$$
i.e. the desired inequality.

Now, for the proof of second inequality of (2.9) we shall apply the inequality
$$L>\sqrt[3]{G\left(\ds\f{A+G}{2}\right)^2},$$
due to the author (see [7]).
This implies
$$\ds\f{1}{2}\left(1+\ds\f{A}{L}\right)<\ds\f{1}{2}\left(1+\sqrt[3]{\ds\f{4A^3}{G(A+G)^2}}\right)=M.$$

By letting $x=\ds\f{A}{G}>1$ in (2.12) of Lemma 2.1 we can deduce
$$M<\sqrt[3]{\ds\f{A}{G}},$$
so
$$\ds\f{1}{2}\left(1+\ds\f{A}{L}\right)<\sqrt[3]{\ds\f{A}{G}}.$$

This may be written also as
$$\ds\f{A+L}{2LA}<\ds\f{1}{\sqrt[3]{A^2G}},$$
 so
$$H(L,A)>\sqrt[3]{A^2G}$$
follows.

{\bf Remark 2.2.}
Therefore, in view of (2.5), (2.6), (2.7) and (2.8), the following
more accurate relations may be written:
$$\ds\f{\sin x}{x}>\ds\f{2L^*}{L^*+1}>\sqrt[3]{\cos x}<\ds\f{3\cos x}{2\cos x+1},\q
x\in \left(0,\ds\f{\pi }{2}\right),
\eqno(2.14)$$
$$\ds\f{\sinh x}{x}>\sqrt[3]{\cosh x}>\ds\f{2P^*}{P^*+1}>\ds\f{3\cosh x}{2\cosh x+1},\q
x>0,
\eqno(2.15)$$
where, as in (2.7) and (2.8),
$$L^*=L(1+\sin x,1-\sin x)
\mbox{ and }
P^*=P(e^x,e^{-x}).$$


\newpage
\setcounter{section}{0}

\bc
{\Large\bf 6. New Refinements of the Huygens and Wilker Hyperbolic and Trigonometric Inequalities}
\ec

\noindent
{\bf AMS Subject Classification:}
26D05, 26D07, 26D99.

\noindent
{\bf Keywords and phrases:}
Inequalities, trigonometric functions, hyperbolic
functions, means and their inequalities.

\section{Introduction}

The famous Huygens, resp. Wilker trigonometric inequalities can be stated as follows:
For any $x\in \left(0,\ds\f{\pi }{2}\right)$ one has
$$2\sin x+\tan x>3x,
\eqno(1.1)$$
resp.
$$\left(\ds\f{\sin x}{x}\right)^2+\ds\f{\tan x}{x}>2.
\eqno(1.2)$$

Their hyperbolic variants are:
For any $x > 0$ hold
$$2\sinh x+\tanh x>3x,
\eqno(1.3)$$
$$\left(\ds\f{\sinh x}{x}\right)^2+\ds\f{\tanh x}{x}>2.
\eqno(1.4)$$

Clearly, (1.4) hold also for $x > 0$, so it holds for any
$x\ne 0$.
In what follows, we shall assume in all inequalities $x>0$
(or $0<x<\ds\f{\pi }{2}$ in trigonometric inequalities).

For references, connections, extensions and history of these inequalities we
quote the recent paper [5].

In what follows, by using the theory of means of two arguments, we will offer
refinements of (1.1) or (1.3), as well as (1.2) or (1.4).

\section{Means of two arguments}

Let $a$, $b$ be positive real numbers.
The logarithmic mean and the identric mean of a and b are defined by
$$L=L(a,b)=\ds\f{a-b}{\ln a-\ln b},
\mbox{ for } a\ne b,\q L(a,a)=a
\eqno(2.1)$$
and
$$I=I(a,b)=\ds\f{1}{e}(a^a/b^b)^{1/(a-b)},\mbox{ for }
a\ne b,\q
I(a,a)=a.
\eqno(2.2.)$$

The Seiffert mean $P$ is defined by
$$P=P(a,b)=\ds\f{a-b}{2\arcsin \left(\ds\f{a-b}{a+b}\right)},
\mbox{ for } a\ne b,\q
P(a,a)=a.
\eqno(2.3)$$

Let
$A=A(a,b)=\ds\f{a+b}{2}$,
$G=G(a,b)=\sqrt{ab}$,
$H=H(a,b)=\ds\f{2ab}{a+b}$
denote the classical arithmetic, geometric, and harmonic means of $a$ and $b$.

There exist many inequalities related to these means.
We quote e.g. [1] for $L$ and $I$, while [2] for the mean $P$.
Recently, E. Neuman and the author [3] have shown that all these means
are the particular cases of the so-called "Schwab-Borchardt mean" (see also [4]).

In what follows, we will use two inequalities which appear in [3] (see Corollary
3.2 of that paper), namely:
$$L>H(P,G),
\eqno(2.4)$$
and
$$P>H(L,A).
\eqno(2.5)$$

Our method will be based on the remark that
$$\ds\f{\sinh x}{x}=L(e^x,e^{-x}),\q x\ne 0
\eqno(2.6)$$
which follows by (2.1), as well as
$$\ds\f{\sin x}{x}=P(1+\cos x,1-\cos x),\q x\in \left(0,\ds\f{\pi }{2}\right)
\eqno(2.7)$$
which may be obtained by definition (2.3).

Note also that
$$G(e^x,e^{-x})=1,\q
A(e^x,e^{-x})=\cosh x
\eqno(2.8)$$
and
$$G(1+\sin x,1-\sin x)=\cos x,\q
A(1+\sin x,1-\sin x)=1.
\eqno(2.9)$$

\section{Main result}
{\bf Theorem 3.1.}
{\it Define
$P^*=P(e^x,e^{-x})$,
$X^*=\ds\f{2\sinh x}{P^*}$.
Then for any $x>0$ one has
$$\tanh x+x>X^*>4x-2\sinh x
\eqno(3.1)$$
and
$$\left(\ds\f{\sinh x}{x}\right)^2+\ds\f{\tanh x}{x}>k^2-2k+3,
\eqno(3.2)$$
where $k=\ds\f{X^*}{2x}$.
}

{\bf Proof.}
Writing inequality (2.4) for
$L=L(e^x,e^{-x})$, etc., and using (2.6) and (2.8), we get:
$$\ds\f{\sinh x}{x}>\ds\f{2P^*}{P^*+1}.
\eqno(3.3)$$

Similarly, for (2.5), we get:
$$P^*>\ds\f{2\sinh x\cosh x}{\sinh x+x\cosh x}.
\eqno(3.4)$$

By (3.4) and (3.3) we can write that
$$\tanh x+x=\ds\f{\sinh x+x\cosh x}{\cosh x}>\ds\f{2\sinh x}{P^*}>4x-2\sinh x,$$
so (3.1) follows.

Now,
$\ds\f{X^*}{x}=2\cdot \ds\f{\sinh x}{x\cdot P^*}=2\cdot \ds\f{L}{P}(e^x,e^{-x})<2$
by the known inequality (see e.g. [2])
$$L<P.
\eqno(3.5)$$

By the right side of (3.1) one has
$$\ds\f{\sinh x}{x}>2-\ds\f{X^*}{2x}=2-k>0.
\eqno(3.6)$$

From the left side of (3.1) we get
$$\ds\f{\tanh x}{x}>\ds\f{X^*}{x}-1=2k-1.
\eqno(3.7)$$

Thus, by (3.6) and (3.7) we can write
$$\left(\ds\f{\sinh x}{x}\right)^2+\ds\f{\tanh x}{x}>(2-k)^2+2k-1
=k^2-2k+3>2$$
by $(k-1)^2>0$.
In fact, $k=\ds\f{L}{P}(e^x,e^{-x})<1$,
so there is strict inequality also in the last inequality.

In case of trigonometric functions, the Huygens inequality is refined in the
same manner, but in case of Wilker's inequality the things are slightly distinct.
\hfill $\square$

{\bf Theorem 3.2.}
{\it Define
$L^*=L(1+\sin x,1-\sin x)$
and
$X^{**}=\ds\f{2\sin x}{L^*}$
for $x\in \left(0,\ds\f{\pi }{2}\right)$.

Then one has the inequalities
$$\tan x+x>X^{**}>4x-2\sin x
\eqno(3.8)$$
and
$$\left(\ds\f{\sin x}{x}\right)^2+\ds\f{\tan x}{x}>
\left\{\ba{lll}
(k^*)^2-2k^*+3, & \mbox{if} & k^*\le 2\medskip \\
2k^*-1, & \mbox{if} & k^*>2
\ea\right.>2
\eqno(3.9)$$
where $k^*=\ds\f{X^{**}}{2x}$.
}

{\bf Proof.}
Applying inequality (2.4) for
$L=L(1+\sin x,1-\sin x)=L^*$,
by (2.7) and (2.9) we get
$$L^*>\ds\f{2\sin x\cos x}{\sin x+x\cos x}.
\eqno(3.10)$$

From (2.5) we get
$$\ds\f{\sin x}{x}>\ds\f{2L^*}{L^*+1}.
\eqno(3.11)$$

Thus,
$$\tan x+x=\ds\f{\sin x+\cos x}{\cos x}>\ds\f{2\sin x}{L^*}>4x-2\sin x$$
by (3.10) and (3.11) respectively.
This gives relation (3.8), which refines (1.1).

Now, by the right side of (3.8) we can write
$$\ds\f{\sin x}{x}>2-\ds\f{\sin x}{x\cdot L^*}=2-k^*.$$

Since
$k^*=\ds\f{P}{L}(1+\sin x,1-\sin x)>1$,
and the upper bound of $k^*$ is $+\infty $ as
$x\to \ds\f{\pi }{2}$, clearly $2-k^*>0$
is not true.

i) If $2-k^*\ge 0$,
then as in the proof of Theorem 3.1, we get from~(3.8)
$$\left(\ds\f{\sin x}{x}\right)^2+\ds\f{\tan x}{x}>(2-k^*)^2+2k^*-1
=(k^*)^2-2k^*+3>2.$$

ii) If $2-k^*<0$, then we use only
$\left(\ds\f{\sin x}{x}\right)^2>0$.
Thus
$$\left(\ds\f{\sin x}{x}\right)^2+\ds\f{\tan x}{x}>\ds\f{\tan x}{x}>2k^*-1.$$

In this case $2k^*-1>3>2$.

Thus, in any case, inequality (3.9) holds true, so a refinement of the
trigonometric Wilker inequality (1.2) is valid.
\hfill $\square$


\newpage
\setcounter{section}{0}

\bc
{\Large\bf 7. On Some New Wilker and Huygens Trigonometric-Hyperbolic Inequalities}
\ec

\begin{abstract}
The hyperbolic counterparts of the Wilker and Huygens trigonometric inequalities
have been introduced by L. Zhu \cite{53-4} and E. Neuman - J. S\'andor \cite{53-2}.
Here we shall study certain new inequalities of Wilker and Huygens type,
involving the trigonometric function $\sin x$ and the hyperbolic function $\sinh x$.
Multiplicative analogues of the stated inequalities are pointed out, too.
\end{abstract}

\noindent
{\bf AMS Subject Classification (2010):} 26D05, 26D07.

\noindent
{\bf Keywords and phrases:}
Inequalities, trigonometric functions, hyperbolic functions.

\section{Introduction}

The famous Wilker inequality for trigonometric functions asserts that
$$\left(\ds\f{\sin x}{x}\right)^2+\ds\f{\tan x}{x}>2,\q 0<x<\ds\f{\pi }{2},
\eqno(1)$$
while the Huygens inequality states that
$$\ds\f{2\sin x}{x}+\ds\f{\tan x}{x}>3,\q 0<x<\ds\f{\pi }{2}.
\eqno(2)$$
For many references related to both inequalities (1) and (2), see e.g. \cite{53-2}.

Recently, S. Wu and H. Srivastava \cite{53-3} have established another inequality of Wilker type, namely
$$\left(\ds\f{x}{\sin x}\right)^2+\ds\f{x}{\tan x}>2,\q 0<x<\ds\f{\pi }{2}
\eqno(3)$$
called in \cite{53-2} as the second Wilker inequality.

L. Zhu \cite{53-4} has proved a hyperbolic version of the first Wilker inequality
$$\left(\ds\f{\sinh x}{x}\right)^2+\ds\f{\tanh x}{x}>2,\q x\ne 0.
\eqno(4)$$

The Huygens trigonometric inequality says that
$$\ds\f{2\sin x}{x}+\ds\f{\tan x}{x}>3,\q 0<x<\ds\f{\pi }{2}.
\eqno(5)$$

For many connections between these inequalities, see \cite{53-2}.

E. Neuman and J. S\'andor \cite{53-2} have recently proved the hyperbolic version
$$\ds\f{2\sinh x}{x}+\ds\f{\tanh x}{x}>3,\q x\ne 0.
\eqno(7)$$

The Cusa-Huygens inequality may be written as
$$\ds\f{2x}{\sin x}+\ds\f{x}{\tan x}>3,\q x\ne 0.
\eqno(6)$$

For example, in \cite{53-2} it is shown that the hyperbolic Huygens inequality (7) implies
the hyperbolic Wilker inequality (4).

In this paper we shall introduce trigonometric-hyperbolic Wilker type inequalities
of the forms
$$\left(\ds\f{\sinh x}{x}\right)^2+\ds\f{\sin x}{x}>2
\eqno(8)$$
or
$$\left(\ds\f{x}{\sin x}\right)^2+\ds\f{x}{\sinh x}>2
\eqno(9)$$
and trigonometric-hyperbolic Huygens type inequalities
$$\ds\f{2\sinh x}{x}+\ds\f{\sin x}{x}>3,
\eqno(10)$$
$$\ds\f{2x}{\sin x}+\ds\f{x}{\sinh x}>3,
\eqno(11)$$
where in all cases $x\in \left(0,\ds\f{\pi }{2}\right]$.

Generalizations, as well as connections will be pointed out, too.

\section{Main results}

First we need some trigonometric-hyperbolic auxiliary results.

{\bf Lemma 1.}
{\it For any $x\in \left[0,\ds\f{\pi }{2}\right]$
$$\cos x\cdot \cosh x\le 1
\eqno(12)$$
and
$$\sin x\cdot \sinh x\le x^2.
\eqno(13)$$

}

{\bf Proof.}
Let us consider the application
$f_1(x)=\cos x\cdot \cosh x-1.$

It is immediate that
$$f'_1(x)=-\sin x\cosh x+\cos x\sinh x,$$
$$f''_1(x)=-2\sin x\sinh x\le 0.$$
Thus $f'_1(x)\le f'_1(0)=0$,
$f_1(x)\le f_1(0)=0$ for $x\ge 0$.
Thus (12) holds true.

Let us now put
$f_2(x)=x^2-\sin x\sinh x$.
As
$$f'_2(x)=2x-\cos x\sinh x-\sin x\cosh x,$$
$$f''_2(x)=2(1-\cos x\cosh x),$$
by (12) one has $f''_2(x)\ge 0$, implying
$f'_2(x)\ge f'_2(0)$ for $x\ge 0$.
Thus $f_2(x)\ge f_2(0)=0$, and (13) follows from $x\in \left[0,\ds\f{\pi }{2}\right]$.

{\bf Remark 1.}
As an application of (12) one can prove that
$$\sinh x\le \tan x\mbox{ for } x\in \left[0,\ds\f{\pi }{2}\right].$$
Indeed, the application $f_3(x)=\sinh x-\tan x$ has a derivative
$$f'_3(x)=\cosh x-\ds\f{1}{\cos^2 x}\le 0,$$
since
$\cos^2 x\cosh x\le \cos x\cosh x\le 1$
by (12).
Thus $f_3(x)\le f_3(0)=0$.

{\bf Lemma 2.}
{\it For any $x\in \left[0,\ds\f{\pi }{2}\right]$ one has
$$\sin^2 x\cdot \sinh x\le x^3,
\eqno(14)$$
and
$$\sin x\cdot \sinh^2 x\ge x^3.
\eqno(15)$$

}

{\bf Proof.}
Since
$\ds\f{\sin x}{x}<1$, and by (13) we can write:
$$\left(\ds\f{\sin x}{x}\right)^2<\ds\f{\sin x}{x}<\ds\f{x}{\sinh x}\mbox{ for }
x\in \left(0,\ds\f{\pi }{2}\right].$$
This implies inequality (14).

For the proof of (15), see \cite{53-1}.

We now are in a position to prove the following generalized trigonometric-hyperbolic
Wilker type inequalities:

{\bf Theorem 1.}
{\it For any $x\in \left(0,\ds\f{\pi }{2}\right]$ and any $t>0$ one has
$$\left(\ds\f{x}{\sin x}\right)^{2t}+\left(\ds\f{x}{\sinh x}\right)^t>2
\eqno(16)$$
and
$$\left(\ds\f{\sinh x}{x}\right)^{2t}+\left(\ds\f{\sin x}{x}\right)^t>2.
\eqno(17)$$

}

{\bf Proof.}
Apply the arithmetic-geometric inequality
$$a+b\ge 2\sqrt{ab}$$
first for
$a=\left(\ds\f{x}{\sin x}\right)^{2t}$, $b=\left(\ds\f{x}{\sinh x}\right)^t$.
Since
$$ab=\left(\ds\f{x^3}{\sin^2 x\sinh x}\right)^t>1$$
by (14), inequality (16) follows.
Put now
$$a=\left(\ds\f{\sinh x}{x}\right)^{2t},\q
b=\left(\ds\f{\sin x}{x}\right)^t,$$
and taking into account of (15), we obtain in a similar fashion, inequality (17).

For $t=1$ from (16) and (17) we get (8), resp. (9).

{\bf Remark 2.}
For the particular case of $t=1$, the second Wilker inequality (3) is stronger than the Wilker
inequality (16) (i.e., in fact (9)).

One has
$$\left(\ds\f{x}{\sin x}\right)^2+\ds\f{x}{\sinh x}>\left(\ds\f{x}{\sin x}\right)^2+\ds\f{x}{\tan x}>2
\eqno(17')$$
Indeed,
$$\ds\f{x}{\sinh x}>\ds\f{x}{\tan x},\q\mbox{as}\q
\sinh x<\tan x\mbox{ for }
x\in \left(0,\ds\f{\pi }{2}\right),$$
see Remark 1.

We now prove Huygens type trigonometric-hyperbolic inequalities.

{\bf Theorem 2.}
{\it For any $x\in \left(0,\ds\f{\pi }{2}\right]$ and any $t>0$ one has
$$2\left(\ds\f{\sinh x}{x}\right)^t+\left(\ds\f{\sin x}{x}\right)^t>3
\eqno(18)$$
and
$$2\left(\ds\f{x}{\sin x}\right)^t+\left(\ds\f{x}{\sinh x}\right)^t>3.
\eqno(19)$$

}

{\bf Proof.}
Apply the arithmetic geometric inequality
$$2a+b=a+a+b\ge 3\sqrt[3]{a^2b}$$
first for
$a=\left(\ds\f{\sinh x}{x}\right)^t$, $b=\left(\ds\f{\sin x}{x}\right)^t$.
Then (18) follows by inequality (15).

By letting
$a=\left(\ds\f{a}{\sin x}\right)^t$, $b=\left(\ds\f{x}{\sinh x}\right)^t$,
via (14) we get the Huygens type inequality (19).

The following result shows a comparison of the deduced inequalities:

{\bf Theorem 3.}
{\it For any $x\in \left(0,\ds\f{\pi }{2}\right]$ one has
$$\left(\ds\f{x}{\sin x}\right)^{2t}+\left(\ds\f{x}{\sinh x}\right)^t
>\left(\ds\f{\sinh x}{x}\right)^{2t}+\left(\ds\f{\sin x}{x}\right)^t>2,
\eqno(20)$$
$$2\left(\ds\f{x}{\sin x}\right)^t+\left(\ds\f{x}{\sinh x}\right)^t
>2\left(\ds\f{\sinh x}{x}\right)^t+\left(\ds\f{\sin x}{x}\right)^t>3.
\eqno(21)$$

}

{\bf Proof.}
The inequality
$a^{2t}+\ds\f{1}{b^t}>b^{2t}+\ds\f{1}{a^t}$ for $a>b$, $t>0$ follows by the strict monotonicity
of the function
$$l(x)=x^{2t}-\ds\f{1}{x^t},\q x>0.$$

Since $l(a)>l(b)$, the above inequality follows.
Let now
$a=\ds\f{x}{\sin x}$, $b=\ds\f{\sinh x}{x}$.
By Lemma 1, (13), since $a>b$, the first inequality of (20) follows.

A similar argument with the function
$$q(x)=2x^t-\ds\f{1}{x^t},\q x>0$$
gives the first part of inequality (21).
Now, by (17) and (18) the complete inequalities (20) and (21) are proved.

Finally, we point out certain "multiplicative analogues",
as well as, improvements of the Wilker or Huygens type inequalities.

Since the method is applicable to all inequalities (1)-(11) or (16)-(20), we present it explicitly
to (1) and (2), or more generally:
$$\left(\ds\f{\sin x}{x}\right)^{2t}+\left(\ds\f{\tan x}{x}\right)^t>2$$
and
$$2\left(\ds\f{\sin x}{x}\right)^t+\left(\ds\f{\tan x}{x}\right)^t>3,$$
where $x\in \left(0,\ds\f{\pi }{2}\right)$.

{\bf Theorem 4.}
{\it For any $x\in \left(0,\ds\f{\pi }{2}\right)$ one has
$$\left[1+\left(\ds\f{\sin x}{x}\right)^{2t}\right]
\left[1+\left(\ds\f{\tan x}{x}\right)^t\right]>4,
\eqno(22)$$
$$\left[1+\left(\ds\f{\sin x}{x}\right)^t \right]^2 \left[1+\left(\ds\f{\tan x}{x}\right)^t\right]>8,
\eqno(23)$$
$$\left(\ds\f{\sin x}{x}\right)^{2t}+\left(\ds\f{\tan x}{x}\right)^t
>2\sqrt{\left[1+\left(\ds\f{\sin x}{x}\right)^{2t}\right]
\left[1+\left(\ds\f{\tan x}{x}\right)^t \right]}-2>2
\eqno(24)$$
$$2\left(\ds\f{\sin x}{x}\right)^{t}+\left(\ds\f{\tan x}{x}\right)^t
>3\sqrt[3]{\left[1+\left(\ds\f{\sin x}{x}\right)^{t}\right]
\left[1+\left(\ds\f{\tan x}{x}\right)^t \right]}-3>3
\eqno(25)$$

}

{\bf Proof.}
We shall apply the following simple refinement of the arithmetic mean - geometric mean inequality:
$$\ds\f{a_1+a_2+\ldots +a_n}{n}\ge
\sqrt[n]{(1+a_1)\ldots (1+a_n)}-1
\ge \sqrt[n]{a_1a_2\ldots a_n}.
\eqno(26)$$

This is discussed e.g. in [?], but we present here a proof of it:
By the known inequality
$$\sqrt[n]{(b_1+a_1)\ldots (b_n+a_n)}
\ge \sqrt[n]{b_1\ldots b_n}+\sqrt[n]{a_1\ldots a_n}$$
for $a_i,b_i>0$ $(i=\ov{1,n})$, applied for $b_i=1$ $(i=\ov{1,n})$ we get
$$\sqrt[n]{(1+a_1)\ldots (1+a_n)}\ge 1+\sqrt[n]{a_1\ldots a_n},$$
which gives the right side of (26).
On the other hand, the arithmetic mean - geometric mean inequality implies
$$\sqrt[n]{(1+a_1)\ldots (1+a_n)}\le 1+\ds\f{a_1+\ldots +a_n}{n},$$
so the left side of (26) follows, too.

Apply now (26) for $n=2$, $a_1=a$, $a_2=b$ in order to deduce
$$a+b\ge 2\sqrt{(1+a)(1+b)}-2>2\sqrt{ab};
\eqno(27)$$
and then for $n=3$, $a_1=a$, $a_2=a$, $a_3=b$ in order to obtain
$$2a+b\ge 3\sqrt[3]{(1+a)^2(1+b)}-3>3\sqrt[3]{a^2b}.
\eqno(28)$$

By letting in (27)
$a=\left(\ds\f{\sin t}{t}\right)^{2t}$, $b=\left(\ds\f{\tan x}{x}\right)^t$,
and using the known inequality
$\left(\ds\f{\sin x}{x}\right)^3>\cos x$
(see \cite{53-2}), we get (22) and (23).
By letting
$a=\left(\ds\f{\sin x}{x}\right)^t$,
$b=\left(\ds\f{\tan x}{x}\right)^t$
in (28), in the same manner, we obtain (24) and (25).

Clearly, the analogous inequalities for all the remaining Wilker or Huygens type inequalities
can be obtained in the same manner.
We omit further details.

{\bf Remark 3.}
Inequality (22) could be called as a "multiplicative Wilker inequality",
and likewise (23) as a "multiplicative Huygens inequality".


\newpage
\setcounter{section}{0}

\bc
{\Large\bf 8. An Optimal Inequality for Hyperbolic and Trigonometric Functions}
\ec

\begin{abstract}
We determine the best positive constants $p$ and $q$ such that
$$\left(\ds\f{1}{\cosh x}\right)^p<\ds\f{\sin x}{x}<\left(\ds\f{1}{\cosh x}\right)^q$$
as well as $p'$ and $q'$ such that
$$\left(\ds\f{\sinh x}{x}\right)^{p'}<\ds\f{2}{\cos x+1}<\left(\ds\f{\sinh x}{x}\right)^{q'}.$$
\end{abstract}



\section{Introduction}

In recent years, many trigonometric and hyperbolic inequalities have attracted
attention of several researchers.
For example, the Huygens, the Cusa-Huygens,
the Wilker trigonometric equalities, or their hyperbolic variants have been studied by many authors.
We quote only the recent papers [1], [4], where many
further references may be found.
In the recent paper [1] there have been considered inequalities satisfied
by both trigonometric and hyperbolic functions.
For example, it was shown that for all $x\in (0,\pi /2)$ one has
$$\ds\f{x^2}{\sinh^2 x}<\ds\f{\sin x}{x}<\ds\f{x}{\sinh x},
\eqno(1.1)$$
$$\ds\f{1}{\cosh x}<\ds\f{\sin x}{x}<\ds\f{x}{\sinh x},
\eqno(1.2)$$
$$\left(\ds\f{1}{\cosh x}\right)^{1/2}<\ds\f{x}{\sinh x}<\left(\ds\f{1}{\cosh x}\right)^{1/4}
\eqno(1.3)$$
for $0<x<1$.

In the recent paper [5] we have determined the best inequalities of type (1.1).
The aim of this paper is to determine the optimal positive constants $p$ and $q$
such that
$$\ds\f{1}{(\cosh x)^p}<\ds\f{\sin x}{x}<\ds\f{1}{(\cosh x)^q},\q
x\in \left(0,\ds\f{\pi }{2}\right)
\eqno(1.4)$$
which is motivated by the above inequalities (1.1)-(1.3).

\section{Main results}
The following auxiliary results will be needed:

{\bf Lemma 2.1.}
{\it For all $x>0$ one has
$$\ln\cosh x>\ds\f{x}{2}\tanh x.
\eqno(2.1)$$

}

{\bf Proof.}
Let us define
$f_1(x)=\ln\cosh x-\ds\f{x}{2}\tanh x$, $x\ge 0$.

A simple computation gives
$$2\cosh^2 x\cdot f'_1(x)=\sinh x\cdot \cosh x-x>0,$$
as it is well known that
$\sinh x>x$, $\cosh x>1$ for $x>0$.
Thus $f_1$ is strictly increasing function, implying
$f_1(x)\ge f_1(0)=0$ for $x\ge 0$, with equality only for $x=0$.
This gives inequality (2.1).
\hfill $\square$

{\bf Lemma 2.2.}
{\it For all $x\in (0,\pi /2)$ one has
$$\ln\ds\f{x}{\sin x}<\ds\f{\sin x-x\cos x}{2\sin x}.
\eqno(2.2)$$

}

{\bf Proof.}
Let
$f_2(x)=\ds\f{\sin x-x\cos x}{2\sin x}-\ln \ds\f{x}{\sin x}$,
$0<x\le \ds\f{\pi }{2}$.

A simple computation gives
$$2x\sin^2 x\cdot f'_2(x)=x^2+x\cdot \sin x\cdot \cos x-2\sin^2 x>0$$
iff
$$t=\ds\f{\sin x}{x}<\ds\f{\cos x+\sqrt{\cos^2 x+8}}{4},
\eqno(2.3)$$
which follows by resolving the second order inequality
$$2t^2-t\cdot \cos x-1<0.$$

Now, by the known Cusa{Huygens inequality (see e.g. [4])
$$\ds\f{\sin x}{x}<\ds\f{\cos x+2}{3},
\eqno(2.4)$$
(2.3) follows, if we can show that
$$\ds\f{\cos x+2}{3}<\ds\f{\cos x+\sqrt{\cos^2 x+8}}{4}.$$

This inequality however becomes, after certain elementary transformations
$$\cos^2 x+1>2\cos x,\q\mbox{i.e.}\q
(\cos x-1)^2>0.$$

Thus $f'_2(x)>0$ for $x>0$, and this implies
$$f_2(x)>f_2(0_+)=\lim\limits_{x\to 0_+}f_2(x)=0,$$
by the known limit formula
$\lim\limits_{x\to 0_+}\ds\f{\sin x}{x}=1$.

The main result of this paper is contained in the following:

{\bf Theorem 2.1.}
{\it The best positive constant $p$ and $q$ in inequality (1.4) are p
$$p=(\ln \pi /2)/(\ln\cosh \pi /2)\approx 0.49
\q\mbox{and}\q
q=\ds\f{1}{3}=0.33\ldots $$

}

{\bf Proof.}
Let us introduce the application
$$h_1(x)=\ds\f{\ln\ds\f{s}{\sin x}}{\ln\cosh x}=\ds\f{f(x)}{g(x)},\q
x\in \left(0,\ds\f{\pi }{2}\right).
\eqno(2.5)$$

Simple computations give
$$f'(x)=\ds\f{\sin x-x\cos x}{x\sin x},\q
g'(x)=\ds\f{\sinh x}{\cosh x},$$
$$(\ln \cosh x)^2  h'_1(x)
=\ds\f{\sin x-x\cos x}{x\sin x} \ln(\cosh x)
-\tanh x \ln \ds\f{x}{\sin x}.
\eqno(2.6)$$

Now, by inequalities
$\sin x>x\cos x$,
$\ds\f{x}{\sin x}>1$,
$\cosh x>1$
and relations (2.1) of Lemma 2.1, resp. (2.2) of Lemma 2.2, by (2.6) we can immediately deduce
$h'_1(x)>0$ for $x>0$.
This shows that, the function $h$ is strictly increasing, so
$$h_1(0_+)<h_1(x)<h_1\left(\ds\f{\pi }{2}\right)
\mbox{ for any } 0<x<\ds\f{\pi }{2}.
\eqno(2.7)$$

Now, elementary computations give
$$h_1(0_+)=\lim\limits_{x\to 0}h_1(x)=\ds\f{1}{2}
\q\mbox{and}\q
h_1(\pi /2)=\ds\f{\ln \pi /2}{\ln\cosh \pi /2}\approx 0.49\ldots ,$$
with the aid of a computer calculation.
Since inequality (1.4) may be written as
$q<h_1(x)<p$,
by (2.7)
$$q=h_1(0_+)
\q\mbox{and}\q
p=h_1(\pi /2)$$
are the best possible constants.
\hfill $\square$

{\bf Remark 2.1.}
The right side inequality of (1.4), with
$q=\ds\f{1}{3}$
follows also from inequality (1.2) and the inequality
$$\ds\f{\sinh x}{x}>\sqrt[3]{\cosh x}
\eqno(2.8)$$
discovered also by I. Lazarevi\'c (see [3], [4]).
We note that, we have shown recently (see [6]) that (2.8) is equivalent also with an inequality in the theory
of means, due to E. B. Leach and M. C. Sholander [2], namely
$$L>\sqrt[3]{G^2A},
\eqno(2.9)$$
where $L=L(a,b)=(b-a)/(\ln b-\ln a)$ $(a\ne b)$
is the logarithmic mean of $a$ and $b$, while
$$G=G(a,b)=\sqrt{ab},\q
A=A(a,b)=\ds\f{a+b}{2}$$
are the well-known geometric and arithmetic means.
See e.g. [6] for these and related means and their inequalities.

We note that inequality (2.1) of Lemma 2.1 follows also from results of the
theory of means.
Let
$$S=S(a,b)=(a^a\cdot b^b)^{1/(a+b)}$$
be a mean studied e.g. in [7].
Now, consider the inequality
$$S<\ds\f{A^2}{G}
\eqno(2.10)$$
(see for improvements and other related inequalities).

Put now
$a=e^x$, $b=e^{-x}$ in (2.10).
As
$$A=A(a,b)=\cosh x,\q
G=G(a,b)=1
\q\mbox{and}\q
S=S(a,b)=e^{x\tanh x},$$
it is immediate that (2.10) becomes (2.1).
From results in [7] we can deduce the following improvement of (2.1):
$$\ln\cosh x>\ds\f{1}{4}[3(x\coth x-1)+x\tanh x]>\ds\f{x}{2}\tanh x.
\eqno(2.11)$$

{\bf Theorem 2.2.}
{\it The best positive constants $p'q$ and $q'$ such that
$$\left(\ds\f{\sinh x}{x}\right)^{p'}<\ds\f{2}{\cos x+1}
<\left(\ds\f{\sinh x}{x}\right)^{q'}
\eqno(2.12)$$
are $p'=\ds\f{3}{2}=1.5$ and
$q'=\ln 2/\ln[(\sinh (\pi /2))/\pi /2]=1.818\ldots $
}

{\bf Proof.}
Let us introduce the application
$$h_2(x)=\ds\f{\ln 2/(\cos x+1)}{\ln\sinh x/x}=\ds\f{F(x)}{G(x)},\q
x\in \left(0,\ds\f{\pi }{2}\right).
\eqno(2.13)$$

As
$F'(x)=\ds\f{\sin x}{\cos x+1}$
and
$G'(x)=\ds\f{x\cosh x-\sinh x}{x\sinh x}$, and
$$G^2(x)\cdot h'_2(x)=\ds\f{x\cosh x-\sinh x}{x\sinh x}
\left(\ln \ds\f{2}{\cos x+1}\right)-\left(\ln \ds\f{\sinh x}{x}\right)
\ds\f{\sin x}{\cos x+1}.
\eqno(2.14)$$

Similarly to Lemma 2.2 one can prove the following inequality:
$$\ln\ds\f{\sinh x}{x}>\ds\f{1}{2}\cdot \ds\f{x\cosh x-\sinh x}{x\sinh x},\q x>0.
\eqno(2.15)$$

We note that (2.15) follows also from the inequality
$$L^2>G\cdot I,
\eqno(2.16)$$
where $L=L(a,b)$ and $I=I(a,b)$
are the logarithmic, resp. identric means of $a$ and $b$, defined by
$$L=(a-b)/(\ln a-\ln b),\q
I=\ds\f{1}{e}(b^b/a^a)^{1/(a-b)}
\mbox{ for } a\ne b.$$
As
$L(e^x,e^{-x})=\ds\f{\sinh x}{x}$,
$I(e^x,e^{-x})=e^{\tanh x-1}$, $G(e^x,e^{-x})=1$,
by inequality (2.16) (see e.g. [7], [8] for this, and related inequalities) we get relation (2.15).

We now prove that
$$a(x)=\ds\f{x}{2}\cdot \ds\f{\sin x}{\cos x+1}-\ln \ds\f{2}{\cos x+1}>0
\mbox{ for } x\in \left(0,\ds\f{\pi }{2}\right).
\eqno(2.17)$$

An immediate computation gives
$$a'(x)=\ds\f{x-\sin x}{2(\cos x+1)}>0,$$
so as $a(0)=0$,
implying relation (2.17).

Now, by (2.15) and (2.17), and taking into account of (2.14) we get
$h'_2(x)>0$ for $x>0$.
Thus $h_2(x)$
is a strictly increasing function, so
$$p'=h_2(0_+)<h_2(x)<h_2(\pi /2)=q'.
\eqno(2.18)$$

A simple computation based on l'Hospital's rule as well as the familiar limits
$$\lim\limits_{x\to 0_+}\ds\f{\sin x}{x}
=\lim\limits_{x\to 0_+}\ds\f{\sinh x}{x}=1$$
imply
$p'=\ds\f{3}{2}=1.5$.
On the other hand, a computer computation shows that
$$q'=\ds\f{\ln 2}{\ln \left(\ds\f{\sinh \pi /2}{\pi /2}\right)}\approx 1.818\ldots $$
This finishes the proof of Theorem 2.2.

{\bf Remark 2.2.}
As
$\ds\f{\cos x+1}{2}=\cos^2 \ds\f{x}{2}$ and
$\sin x=2\sin\ds\f{x}{2}\cos\ds\f{x}{2}$,
$\sin \ds\f{x}{2}<\ds\f{x}{2}$,
$\tan \ds\f{x}{2}>\ds\f{x}{2}$,
it is easy to verify that
$$\left(\ds\f{\sin x}{x}\right)^2<\ds\f{\cos x+1}{2}<\ds\f{\sin x}{x}.
\eqno(2.19)$$
By inequality (1.1) these imply
$$\ds\f{\sinh x}{x}<\ds\f{2}{\cos x+1}<\left(\ds\f{\sinh x}{x}\right)^4.
\eqno(2.20)$$
The above proved theorem offers the sharpest version of this inequality.


\newpage
\setcounter{section}{0}

\bc
{\Large\bf 9. Two Sharp Inequalities for Trigonometric and Hyperbolic Functions}
\ec

\begin{abstract}
We determine the best positive constants $p$ and $q$ such that
$$(\sinh x/x)^p<x/sin x<(\sinh x/x)^q.$$
\end{abstract}



\section{Introduction}
The trigonometric and hyperbolic inequalities have been in the last years in the
focus of many researchers.
For many results and a long list of references we quote the papers [1], [2].
In paper [1], the following interesting inequalities have been proved:

{\bf Theorem 1.1.}
{\it For any $x\in \left(0,\ds\f{\pi }{2}\right)$,
$$\ds\f{\sinh x}{x}<\ds\f{x}{\sin x}<\left(\ds\f{\sinh x}{x}\right)^2.
\eqno(1.1)$$

}

The aim of this paper is to find the best possible form of inequalities (1.1)
in the sense of determination of greatest $p > 0$ and least $q > 0$ such that
$$\left(\ds\f{\sinh x}{x}\right)^p<\ds\f{x}{\sin x}<\left(\ds\f{\sinh x}{x}\right)^q.
\eqno(1.2)$$

\section{Main results}
First we state two auxiliary results:

{\bf Lemma 2.1.}
{\it For any $x > 0$ one has
$$\ln \ds\f{\sinh x}{x}>\ds\f{x\cosh x-\sinh x}{2\sinh x}.
\eqno(2.1)$$

}

{\bf Proof.}
Let
$$L(a,b)=\ds\f{a-b}{\ln a-\ln b},\q
G(a,b)=\sqrt{ab},\q
I(a,b)=\ds\f{1}{e}(b^b/a^a)^{1/(b-a)}$$
for $a\ne b$ positive numbers.
These are called as the logarithmic, geometric and identric means of $a$ and $b$
(see e.g. [3] for $L$, $I$ and related means).
The inequality
$$L^2>G\cdot I
\eqno(2.2)$$
is due to H. Alzer [4].
Apply (2.2) for
$a=e^x$, $b=e^{-x}$.
Then it is easy to see that
$$L=L(a,b)=\ds\f{\sinh x}{x},\q
G=G(a,b)=1,\q
I=I(a,b)=e^{x\coth x-1}.$$

After these substitutions, inequality (2.1) follows at once.

{\bf Lemma 2.2.}
{\it For any $x\in \left(0,\ds\f{\pi }{2}\right)$ one has
$$\ln \ds\f{x}{\sin x}<\ds\f{\sin x-x\cos x}{2\sin x}.
\eqno(2.3)$$

}

{\bf Proof.}
Let
$a(x)=\ds\f{\sin x-x\cos x}{2\sin x}-\ln \ds\f{x}{\sin x}$,
$0<x<\ds\f{\pi }{2}$.

A simple computation gives
$$a'(x)=\ds\f{x^2+x\sin x\cos x-2\sin^2 x}{2x\sin^2 x}>0,$$
if $2\sin^2 x<x^2+x\sin x\cos x$.

As this may be written also as
$$2t^2-t\cos x-1<0,$$
where $t=\ds\f{\sin x}{x}$,
by resolving this inequality of second variable in $t$, as $t > 0$, this is equivalent with
$$\ds\f{\sin x}{x}<\ds\f{\cos x+\sqrt{\cos^2 x+8}}{4}.
\eqno(2.4)$$

But this follows by the "Cusa-Huygens inequality" (see [2])
$$\ds\f{\sin x}{x}<\ds\f{2+\cos x}{3}
\eqno(2.5)$$
as
$\ds\f{2+\cos x}{3}<\ds\f{\cos x+\sqrt{\cos^2 x+8}}{4}$,
as a simple verification shows.
Therefore, the application $a(x)$ is strictly increasing, and as
$\lim\limits_{x\to 0_+}a(x)=0$,
relation (2.3) follows.
\hfill $\square$

{\bf Remark 2.1.}
Relation (2.1) could have been proved via auxiliary functions, too,
but we opted for a proof with means, as this was the first method of discovery of
(2.1) by the author.

We note that, by considering Seiffert's mean P (see [5]) defined by
$$P(a,b)=\ds\f{a-b}{2\arcsin \left(\ds\f{a-b}{a+b}\right)},\q a\ne b$$
and using the substitution
$a=A(1+\sin x)$, $B=A(1-\sin x)$,
where $A=\ds\f{a+b}{2}$ and
$x=\arcsin \ds\f{a-b}{a+b}$ for $0<a<b$;
inequality (2.3) could be written as an inequality for means:
$$P^2>A\cdot X
\eqno(2.6)$$
where
$X=Ae^{G/P-1}$, where $X=X(a,b)$
a new mean (introduced for the first time here).

Indeed, (2.6) follows by remarking that
$P(a,b)=A\cdot \ds\f{\sin x}{x}$,
etc. and introducing in (2.3).

We state the main result of this paper:

{\bf Theorem 2.1.}
{\it The application
$$h(x)=\ds\f{\ln\ds\f{x}{\sin x}}{\ln \ds\f{\sinh x}{x}},\q
x\in \left(0,\ds\f{\pi }{2}\right)
\eqno(2.7)$$
is strictly increasing.

As a corollary, the best constants in inequality (1.2) are
$$p=1
\q\mbox{and}\q
q=\ds\f{\ln(\pi /2)}{\ln((\sinh \pi /2)/\pi /2)}\approx 1.18\ldots $$

}

{\bf Proof.}
Put
$f(x)=\ln\ds\f{x}{\sin x}$,
$g(x)=\ln \ds\f{\sinh x}{x}$.
Simple computations give
$$f'(x)=\ds\f{\sin x-x\cos x}{x\sin x},\q
g'(x)=\ds\f{x\cosh x-\sinh x}{x\sinh x}$$
and
$$xg^2(x)h'(x)
=\left(\ds\f{\sin x-x\cos x}{\sin x}\right)\ln \ds\f{\sinh x}{x}
-\left(\ds\f{x\cosh x-\sinh x}{x\sinh x}\right)\ln \ds\f{x}{\sin x}.
\eqno(2.8)$$

By (2.1) and (2.3) the right side of (2.8) is strictly positive, implying
$h'(x)>0$.
Thus $h$ is strictly increasing.
This implies
$$\lim\limits_{x\to 0_+}h(x)<h(x)<h\left(\ds\f{\pi }{2}\right)
=\ds\f{\ln \pi /2}{\ln \left(\ds\f{\sin \pi /2}{\pi /2}\right)}.$$

A simple computation shows
$\lim\limits_{x\to 0_+}h(x)=1$,
thus
$$p=1<h(x)<\ds\f{\ln \pi /2}{\ln \ds\f{(\sinh \pi /2)}{\pi /2}}=q.
\eqno(2.9)$$

Clearly, the inequality (1.2) maybe written also as
$$p<h(x)<q,$$
so (2.9) offer the best possible constants $p$ and $q$.

{\bf Remark 2.2.}
A computer computation shows
$q\approx 1.18\ldots <2$,
thus inequality (1.1) is improved on the right-hand side.

\section{An application}
A famous inequality of Wilker (see [2] for many connections with other inequalities) states that
$$\left(\ds\f{\sin x}{x}\right)^2+\ds\f{\tan x}{x}>2,\q
0<x<\ds\f{\pi }{2}.
\eqno(3.1)$$

In what follows, we will obtain a Wilker type inequality for the functions
$\sin x$ and $\sinh x$:

{\bf Theorem 3.1.}
{\it Let
$q=\ds\f{\ln \pi /2}{\ln \ds\f{(\sinh \pi /2)}{\pi /2}}\approx 1.18\ldots $

Then
$$\ds\f{\sin x}{x}+\left(\ds\f{\sinh x}{x}\right)^q>q+1.
\eqno(3.2)$$

}

{\bf Proof.}
Apply the pondered arithmetic-geometric inequality
$$\lambda a+(1-\lambda )b>a^\lambda b^{1-\lambda }$$
for $a>0$, $b>0$, $a\ne b$, $\lambda \in (0,1)$.
Put
$\lambda =\ds\f{1}{q+1}$ and $a=\ds\f{\sin x}{x}$, $b=\ds\f{\sinh x}{x}$.
Then, as the right side of (1.2) may be written as
$ab^\lambda >1$, and as
$$a+qb>(ab^q)^{1/(q+1)}(q+1)>q+1.
\eqno\square $$

The following result will be a Wilker type inequality for the functions
$\ds\f{\sin x}{x}$ and $\ds\f{\sinh x}{x}$:

{\bf Theorem 3.2.}
{\it
$$\left(\ds\f{\sinh x}{x}\right)^q+\ds\f{\sin x}{x}>2
\mbox{ for }
x\in \left(0,\ds\f{\pi }{2}\right),
\eqno(3.3)$$
where $q$ is the same as in Theorem 3.1.
}

{\bf Proof.}
The Bernoulli inequality states that (see e.g. [6])
$$(1+t)^\alpha >1+\alpha t,\mbox{ for } t>0,\ \alpha >0.
\eqno(3.4)$$

Put $t=\ds\f{\sinh x}{x}-1>0$, $\alpha =q$ in (3.4).
Then by this inequality, combined with inequality (2.9) implies (3.3).
\hfill $\square$

{\bf Remark 3.1.}
As
$\left(\ds\f{\sinh x}{x}\right)^2>\left(\ds\f{\sinh x}{x}\right)^q$,
we get from (3.3) a much weaker, but more familiar form of Wilker's inequality:
$$\left(\ds\f{\sinh x}{x}\right)^2+\ds\f{\sin x}{x}>2,\q x\in \left(0,\ds\f{\pi }{2}\right).
\eqno(3.5)$$


\newpage
\setcounter{section}{0}

\bc
{\Large\bf 10. On Two New Means of Two Variables}
\ec

\begin{abstract}
Let $A$, $G$ and $L$ denote the arithmetic, geometric resp. logarithmic
means of two positive number, and let $P$ denote the Seiffert mean.
We study the properties of two new means $X$ resp. $Y$, defined by
$$X=A\cdot e^{G/P-1}
\q\mbox{and}\q
Y=G\cdot e^{L/A-1}.$$
\end{abstract}



\section{Introduction}

Let $a$, $b$ be two positive numbers.
The logarithmic and identric means of $a$ and $b$ b are defined by
$$\ba{l}
L=L(a,b)=\ds\f{a-b}{\ln a-\ln b}\ (a\ne b),\q L(a,a)=a;\medskip \\
I=I(a,b)=\ds\f{1}{e}(b^b/a^a)^{1/(b-a)}\ (a\ne b),\q I(a,a)=a.
\ea
\eqno(1.1)$$

The Seiffert mean $P$ is defined by
$$P=P(a,b)=\ds\f{b-a}{2\arcsin \ds\f{b-a}{a+b}}\ (a\ne b),\q
P(a,a)=a.
\eqno(1.2)$$

Let
$$A=A(a,b)=\ds\f{a+b}{2},\q
G=G(a,b)=\sqrt{ab}\q\mbox{and}\q
H=H(a,b)=\ds\f{2ab}{a+b}$$
denote the classical arithmetic, geometric, resp. harmonic means of $a$ and $b$.
There exist many papers which study properties of these means.
We quote e.g. [1], [2] for the identric and logarithmic means, and [3] for the mean $P$.

The means $L$, $I$ and $P$ are particular cases of the "Schwab-Borchardt mean",
see [4], [5] for details.

The aim of this paper is the study of two new means, which we shall denote
by
$X=X(a,b)$ and $Y=Y(a,b)$,
defined as follows:
$$X=A\cdot e^{\f{G}{P}-1},
\eqno(1.3)$$
resp.
$$Y=G\cdot e^{\f{L}{A}-1}.
\eqno(1.4)$$

Clearly
$X(a,a)=Y(a,a)=a$,
but we will be mainly interested for properties of these means for $a\ne b$.

\section{Main results}
{\bf Lemma 2.1.}
{\it The function
$f(t)=te^{\f{1}{t}-1}$, $t>1$
is strictly increasing.
For all $t>0$, $t\ne 0$ one has $f(t)>1$.
For $0<t<1$, $f$ is strictly decreasing.
As a corollary, for all $t>0$, $t\ne 1$ one has
$$1-\ds\f{1}{t}<\ln t<t-1.
\eqno(2.1)$$

}

{\bf Proof.}
As
$\ln f(t)=\ln t+\ds\f{1}{t}-1$,
we get
$$\ds\f{f'(t)}{f(t)}=\ds\f{t-1}{t^2},$$
so $t_0=1$
will be a minimum point of $f(t)$, implying
$f(t)\ge f(1)=1$, for any $t>0$,
with equality only for $t = 1$.
By taking logarithm, the left side of (2.1) follows.
Putting $1/t$ in place of $t$, the left side of (2.1) implies the right side inequality.
\hfill $\square$

{\bf Theorem 2.1.}
{\it For $a\ne b$ one has
$$G<\ds\f{A\cdot G}{P}<X<\ds\f{A\cdot P}{2P-G}<P.
\eqno(2.2)$$

}

{\bf Proof.}
Applying (2.1) for $t=\ds\f{X}{A}$ ($\ne 1$, as $G\ne P$ for $a\ne b$),
and by taking into account of (1.3) we get the middle inequalities of (2.2).
As it is well known that (see [3])
$$\ds\f{A+G}{2}<P<A,
\eqno(2.3)$$
the first inequality of (2.3) implies the last one of (2.2), while the second inequality
of (2.3) implies the first one of (2.2).

In a similar manner, the following is true:

{\bf Theorem 2.2.}
{\it For $a\ne b$ one has
$$H<\ds\f{L\cdot G}{A}<Y<\ds\f{G\cdot A}{2A-L}<G.
\eqno(2.4)$$

}

{\bf Proof.}
Applying (2.1) for $t=\ds\f{Y}{G}$
by (1.4) we can deduce the second and third inequalities of (2.4).
Since
$H=\ds\f{G^2}{A}$,
the first and last inequality of (2.4) follows by the known inequalities (see e.g. [1] for references)
$$G<L<A.
\eqno(2.5)$$
\hfill $\square$

The second inequality of (2.2) can be strongly improved, as follows:

{\bf Theorem 2.3.}
{\it For $a\ne b$ one has
$$1<\ds\f{L^2}{G\cdot I}<\ds\f{L}{G}\cdot e^{\f{G}{L}-1}<\ds\f{X\cdot P}{A\cdot G}.
\eqno(2.6)$$

}

{\bf Proof.}
As $L < P$ (due to Seiffert; see [3] for references) and
$$f\left(\ds\f{P}{G}\right)=X\cdot \ds\f{P}{A\cdot G},
\eqno(2.7)$$
where $f$ is defined in Lemma 2.1, and by taking into account of the inequality
(see [1])
$$\ds\f{L}{I}<e^{\f{G}{L}-1},
\eqno(2.8)$$
by the monotonicity of $f$ one has
$$f\left(\ds\f{P}{G}\right)=\ds\f{L}{G}\cdot e^{\f{G}{L}-1}
>\ds\f{L}{G}\cdot \ds\f{L}{I}=\ds\f{L^2}{G\cdot I}.
\eqno(2.9)$$

By an inequality of Alzer (see [1] for references) one has
$$L^2>G\cdot I,
\eqno(2.10)$$
thus all inequalities of (2.6) are established.

The following estimates improve the left side of (2.4):

{\bf Theorem 2.4.}
{\it For $a\ne b$,
$$H<\ds\f{G^2}{I}<\ds\f{L\cdot G}{A}<\ds\f{G\cdot (A+L)}{3A-L}<Y.
\eqno(2.11)$$

}

{\bf Proof.}
Since $H=\ds\f{G^2}{A}$,
the first inequality of (2.11) follows by the known inequality $I < A$
(see [1] for references).
The second inequality of (2.11) follows by
another known result of Alzer (see [1] for references, and [3] for improvements)
$$A\cdot G<L\cdot I.
\eqno(2.12)$$

Finally, to prove the last inequality of $Y$, remark that the logarithmic mean of $Y$ and $G$ is
$$L(Y,G)=\ds\f{Y-G}{\ln Y/G}=\ds\f{(G-Y)A}{A-L}.
\eqno(2.13)$$

Now, by the right side of (2.5) applied to $a = Y$, $b = G$ we have
$$L(Y,G)<(Y+G)/2,$$
so
$$2A(G-Y)<(A-L)(Y+G),$$
which after some transformations gives the desired inequality.
\hfill $\square$

Similarly to (2.11) we can state:

{\bf Theorem 2.5.}
$$\ds\f{A\cdot G}{P}<\ds\f{A(P+G)}{3P-G}<X.
\eqno(2.14)$$

{\bf Proof.}
$L(X,A)=(X-A)\log X/A=\ds\f{(A-X)P}{P-G}<\ds\f{X+A}{2}$,
so after simple computations we get the second inequality of (2.14).
The first inequality becomes
$$(P-G)^2>0.
\eqno\square $$

A connection between the two means $X$ and $Y$ is provided by:

{\bf Theorem 2.6.}
{\it For $a\ne b$,
$$A^2\cdot Y<P\cdot L\cdot X.
\eqno(2.15)$$

}

{\bf Proof.}
By using the inequality (see [3])
$$\ds\f{A}{L}<\ds\f{P}{G},
\eqno(2.16)$$
and remarking that
$$f\left(\ds\f{A}{L}\right)=\ds\f{A}{L\cdot G}\cdot Y,
\eqno(2.17)$$
by the monotonicity of $f$ one has
$f\left(\ds\f{A}{L}\right)<f\left(\ds\f{P}{G}\right)$,
so by (2.7) and (2.17) we can deduce inequality (2.15).
\hfill $\square$

{\bf Remark 2.1.}
By the known identity (see [1], [2])
$$\ds\f{I}{G}=e^{\f{A}{L}-1}
\eqno(2.18)$$
and the above methods one can deduce the following inequalities (for $a\ne b$):
$$1<\ds\f{L\cdot I}{A\cdot G}<\ds\f{G}{P}\cdot e^{\f{P}{G}-1}.
\eqno(2.19)$$

Indeed, as
$f\left(\ds\f{L}{A}\right)=\ds\f{L}{A}\cdot e^{\f{A}{L}-1}=\ds\f{L\cdot I}{A\cdot G}>1$
we reobtain inequality (2.12).
On the other hand, by (2.16) we can write, as
$1>\ds\f{L}{A}>\ds\f{G}{P}$
that
$f\left(\ds\f{L}{A}\right)<f\left(\ds\f{G}{P}\right)$
i.e. the complete inequality (2.19) is established.

{\bf Theorem 2.7.}
{\it For $a\ne b$
$$X<A\left[\ds\f{1}{e}+\left(1-\ds\f{1}{e}\right)\ds\f{G}{P}\right]
\eqno(2.20)$$
and
$$Y<G\left[\ds\f{1}{e}+\left(1-\ds\f{1}{e}\right)\ds\f{L}{A}\right].
\eqno(2.21)$$

}

{\bf Proof.}
The following auxiliary result will be used:

{\bf Lemma 2.2.}
{\it For the function $f$ of Lemma 2.1, for any $t > 1$ one has
$$f(t)<\ds\f{1}{e}(t+e-1)
\eqno(2.22)$$
and
$$f(t)<\ds\f{1}{e}\left(t+\ds\f{1}{2t}+e-\ds\f{3}{2}\right)<\ds\f{1}{e}(t+e-1).
\eqno(2.23)$$

}

{\bf Proof.}
By the series expansion of $e^x$ and by $t > 1$, we have
$$f(t)=\ds\f{1}{e}\left(t+1+\ds\f{1}{2t}+\ds\f{1}{3!t^2}+\ds\f{1}{4!t^3}+\ldots \right)$$
$$=\ds\f{1}{e}\left(t+\ds\f{1}{1!}+\ds\f{1}{2!}+\ldots \right)=\ds\f{1}{e}(t+e-1),$$
so (2.22) follows.

Similarly,
$$f(t)=\ds\f{1}{e}\left(t+1+\ds\f{1}{2t}+\ds\f{1}{3!t^2}+\ds\f{1}{4! t^3}+\ldots \right)$$
$$<\ds\f{1}{e}\left[t+1+\ds\f{1}{2t}+e-\left(1+\ds\f{1}{1!}+\ds\f{1}{2!}\right)\right]
=\ds\f{1}{e}\left(t+\ds\f{1}{2t}+e-\ds\f{3}{2}\right),$$
so (2.23) follows as well.

Now, (2.20) follows by (2.22) and (2.7), while (2.21) by (2.23) and (2.17).
\hfill $\square$

{\bf Theorem 2.8.}
{\it For $a\ne b$ one has
$$P^2>A\cdot X
\eqno(2.24)$$

}

{\bf Proof.}
Let $x\in \left(0,\ds\f{\pi }{2}\right)$.
In the recent paper [7] we have proved the following trigonometric inequality:
$$\ln \ds\f{x}{\sin x}<\ds\f{\sin x-\cos x}{2\sin x}.
\eqno(2.25)$$

Remark that by (1.2) one has
$$P(1+\sin x,1-\sin x)=\ds\f{\sin x}{x},$$
$$A(1+\sin x,1-\sin x)=1,\q
G(1+\sin x,1-\sin x)=\cos x,$$
so (2.24) may be rewritten also as
$$P^2(1+\sin x,1-\sin x)>A(1+\sin x,1-\sin x)\cdot X(1+\sin x,1-\sin x).
\eqno(2.26)$$

For any $a,b>0$ one can find $x\in \left(0,\ds\f{\pi }{2}\right)$ and $k>0$ such that
$$a=(1+\sin x)k,\q
b=(1-\sin x)k.$$

Indeed, let
$k=\ds\f{a+b}{2}$ and
$x=\arcsin \ds\f{a-b}{a+b}$.

Since the means $P$, $A$ and $X$
are homogeneous of order one, by multiplying (2.26) by $k$, we get the general
inequality (2.24).
\hfill $\square$

{\bf Corollary 2.1.}
$$P^3>\ds\f{A^2L^2}{I}>A^2G.
\eqno(2.27)$$

{\bf Proof.}
By (2.6) of Theorem 2.3 and (2.24) one has
$$\ds\f{L^2A}{I\cdot P}<X<\ds\f{P^2}{A},
\eqno(2.28)$$
so we get
$P^3>\ds\f{A^2L^2}{I}>A^2G$ by inequality (2.10).
\hfill $\square$

{\bf Remark 2.2.}
Inequality (2.27) offers an improvement of
$$P^3>A^2G
\eqno(2.29)$$
from paper [3].
We note that further improvements, in terms of $A$ and $G$ can
be deduced by the "sequential method" of [3].
For any application of (2.29),
put
$a=1+\sin x$, $b=1-\sin x$ in (2.29) to deduce
$$\ds\f{\sin x}{x}>\sqrt[3]{\cos x},\q
x\in \left(0,\ds\f{\pi }{2}\right),
\eqno(2.30)$$
which is called also the Mitrinovi\'c-Adamovi\'c inequality (see [6]).

Since (see [3])
$$P<\ds\f{2A+G}{3},
\eqno(2.31)$$
by the above method we can deduce
$$\ds\f{\sin x}{x}<\ds\f{\cos x+2}{3},
\eqno(2.32)$$
called also as the Cusa-Huygens inequality.
For details on such trigonometric or related hyperbolic inequalities, see [6].

\end{document}